\def\mapright#1{\smash{\mathop{\longrightarrow}\limits^{#1}}}

\def\mapup#1{\uparrow\rlap{$\vcenter{\hbox{$\scriptstyle#1$}}$}}
\def\mapupleft#1{\uparrow\llap
                   {$\vcenter{\hbox{$\scriptstyle#1$}}\;\;$}}

\documentclass{amsart}

\usepackage{amscd, amsmath, amsthm, amssymb}

\newtheorem{theorem}{Theorem}[section]
\newtheorem{lemma}[theorem]{Lemma}
\newtheorem{proposition}[theorem]{Proposition}

\newtheorem{corollary}[theorem]{Corollary}

\theoremstyle{definition}     
\newtheorem{definition}[theorem]{Definition}

\newtheorem{example}[theorem]{Example}

\newtheorem{claim}[theorem]{Claim}

\theoremstyle{remark}
\newtheorem{remark}[theorem]{Remark}

\numberwithin{equation}{section}

\begin{document}

\title[Fermat quartic surface]{A characterization of the 
Fermat quartic K3 surface by means of finite symmetries}

\author[K. Oguiso]{Keiji Oguiso}

\address{
Department of Mathematical Sciences,
University of Tokyo, Komaba Meguro-ku,
Tokyo 153-8914, Japan
}
\email{oguiso@ms.u-tokyo.ac.jp}


\begin{abstract}
We shall characterize the Fermat quartic K3 surface, among all K3 surfaces,  by means of 
its finite group symmetries. 
\end{abstract}

\maketitle


\setcounter{section}{0}
\section{ Introduction}

The aim of this note is to characterize the Fermat quartic surface, among all 
complex K3 surfaces, in terms of finite group symmetries. Our main result is 
theorem 1.2.

Throughout this note, we shall work over the complex number 
field $\mathbf C$. By a K3 surface, we mean a simply-connected 
smooth complex surface $X$ which admits a nowhere vanishing global 
holomorphic 2-from 
$\omega_{X}$. As it is well known, K3 surfaces form $20$-dimensional 
families and projective ones 
form countably many $19$-dimensional families [PSS]. Among them, one of the 
simplest examples is the Fermat quartic surface: 
$$\iota : X_{4} := (x_{1}^{4} + x_{2}^{4} + x_{3}^{4} + x_{4}^{4} = 0) 
\subset \mathbf P^{3}\,\, .$$

From the explicit form, we see that $X_{4}$ admits a fairly large 
projective transformation group, namely, 

$$\tilde{F}_{384} := (\mu_{4}^{4} : S_{4})/\mu_{4} = 
(\mu_{4}^{4}/\mu_{4}) : S_{4}\,\, .$$ 
Here the symbol $A:B$ means a semi-direct product($A$ being normal) 
and $\mu_{I} 
:= \langle \zeta_{I} \rangle$ ($\zeta_{I}:= 
e^{2\pi i/I}$), the multiplicative subgroup of order $I$ of 
$\mathbf C^{\times}$. 
This group $\tilde{F}_{384}$ is a solvable group of order 
$4^{3}\cdot 4! = 2^{9}\cdot 3$. The action by $\tilde{F}_{384}$ 
is an ovbious one, that is, 
$\mu_{4}^{4}$, or 
$\mu_{4}^{4}/\mu_{4}$ acts on $X_{4}$ diagonally 
and $S_{4}$ acts as the permutation of the coordinates.

Let $\tilde{F}_{128}$ be a Sylow $2$-subgroup of $\tilde{F}_{384}$. 
Then $\tilde{F}_{128}$ 
is a nilpotent group of order $2^{9}$. We have an action by 
$\tilde{F}_{128}$ on $X_{4}$ which is a restriction of the action 
by $\tilde{F}_{384}$.  
We call the action 
$$\iota_{384} : \tilde{F}_{384} \times X_{4} \longrightarrow X_{4}\,\, ,\,\, 
\text{resp.}\,\,
\iota_{128} : \tilde{F}_{128} \times X_{4} \longrightarrow X_{4}$$ 
defined here, the standard action by $\tilde{F}_{384}$ (resp. 
by $\tilde{F}_{128}$) on $X_{4}$. 
Note that 
$\tilde{F}_{384}$ has exactly $3$ Sylow $2$-subgroups, corresponding to the $3$  Sylow $2$-subgroups ($\simeq D_{8}$) of $S_{4}$. However, 
they are conjugate one 
another by the Sylow theorem, and 
their standard actions on $X_{4}$ are isomorphic to one another 
in the sense 
below. The group $\tilde{F}_{128}$ is also interesting from 
the view of Mukai's 
classification of 
symplectic K3 groups [Mu]. In fact, 
it is an  
extension of a Sylow $2$-subgroup $F_{128}$ of the Mathieu 
group $M_{23}$ by $\mu_{4}$ from the right (see also section 2). 

\begin{definition} \label{def:K3gp}
We call a finite group $G$ a K3 group (on $X$) if there is a faithful 
action by $G$ 
on $X$, say, $\rho : G \times X \longrightarrow X$. 
Let $G_{i}$ be a K3 group on $X_{i}$ acting by $\rho_{i} : G_{i} 
\times X_{i} \longrightarrow X_{i}$ 
($i = 1, 2$). We say that $(G_{i}, X_{i}, \rho_{i})$ are isomorphic if 
there are a group isomorphism $f : G_{1} \simeq G_{2}$ and an isomorphism 
$\varphi : X_{1} \simeq X_{2}$ such that the following diagram commutes:
$$
\begin{matrix}
G_{2} \times X_{2} & \mapright{\rho_{2}} & X_{2} \cr
\mapupleft{f \times \varphi \; } &  & \mapup{\varphi} \cr
G_{1} \times X_{1} & \mapright{\rho_{1}} & X_{1} \cr
\end{matrix}
$$
\end{definition}
The aim of this note is to show the following theorem:
\begin{theorem} \label{thm:mainThm}
\begin{list}{}{
\setlength{\leftmargin}{10pt}
\setlength{\labelwidth}{6pt}
}
\item[(1)] Let $G$ be a solvable K3 group on $X$ acting by 
$\rho : G \times X \longrightarrow X$. Then $\vert G \vert \le 2^{9}\cdot 3$. 
Moreover, if 
$\vert G \vert = 2^{9} \cdot 3\, (= 1536)$, then $\text{Pic}(X)^{G} = \mathbf Z H$, $(H^{2}) = 4$ 
and $(G, X, \rho) \simeq (\tilde{F}_{384}, X_{4}, \iota_{384})$, 
the standard action by $\tilde{F}_{384}$ on the Fermat 
quartic surface $X_{4}$.  
\item[(2)] Let $G$ be a nilpotent K3 group on $X$ acting by 
$\rho : G \times X \longrightarrow X$. Then 
$\vert G \vert \le 2^{9}$. Moreover, if 
$\vert G \vert = 2^{9}$, then $\text{Pic}(X)^{G} = \mathbf Z H$, 
$(H^{2}) = 4$ and 
$(G, X, \rho) \simeq (\tilde{F}_{128}, X_{4}, \iota_{128})$, 
the standard action by $\tilde{F}_{128}$ on $X_{4}$.  
\end{list}
\end{theorem}

The most basic class of finite groups is the class of cyclic 
groups of prime order. 
This class is extended to the following sequences of important 
classes of groups of rather different nature: 
$$(\text{abelian\, groups}) \subset (\text{nilpotent\, groups})
 \subset (\text{solvable\, groups})\,\, ;$$ 
$$(\text{quasi-simple\, non-commutative\, groups}) \subset 
(\text{quasi-perfect\, groups})\,\, ,$$
where a quasi-simple non-commutative group (resp. a quasi-perfect group) 
is a group which is an 
extension of a simple non-commutative group (resp. a perfect group) by a 
cyclic group from the right.

From the view of these sequences, our theorem is regarded as both an analogy 
and a counter part 
of previous works of 
Kondo [Ko2] for the quasi-perfect 
K3 group $M_{20}:\mu_{4}$, which is also the group of the maximum order among 
K3 groups, and Zhang and the author [OZ] for the quasi-simple non-commutative 
group $L_{2}(7) \times \mu_{4}$. In terms of the coarse moduli space 
$\mathcal M_{4}$ of 
quasi-polarized K3 surfaces of 
degree $4$, our theorem says that the large stabilizer subgroups 
$\tilde{F}_{384}$, 
$\tilde{F}_{128}$ identify the point corresponding 
to the Fermat K3 surface (naturally polarized by $\iota$) 
in $\mathcal M_{4}$. 
However, our theorem claims much more, because 
we do not assume in a priori a degree of invariant polarization. 
Indeed, as 
in [OZ], 
the determination of the degree of invariant polarization is one of 
the key steps in our proof (proposition (3.3) and section 6). For 
this step, we apply Kondo's embedding theorem [Ko1, Lemmas 5, 6] (see also 
(6.2) and (6.3)) as in [Ko2] and [OZ]. 
On the 
other hand, our 
theorem can be also viewed as a characterization of a 
$2$-group $\tilde{F}_{128}$ 
by means of geometry. It might be worth noticing the following table in [Ha, Page 11] of 
the number 
$p(n)$ of isomorphism 
classes of $2$-groups of order $2^{n}$: 

$$
\begin{tabular}{*{10}{|c}{|}} 
\hline 
$n$ & 1 & 2 & 3 & 4 & 5 & 6 & 7 & 8 & 9 \\ 
\hline 
$p(n)$ & 1 & 2 & 5 & 14 & 51 & 267 & 2328 & 56092 & 10494213 \\
\hline 
\end{tabular} 
$$

Section 2 is a summary of known results, relevant to us, about K3 groups 
from Nikulin [Ni] and Mukai [Mu]. In section 3, we reduce our theorem to 
three propositions (propositions (3.1), (3.2), and (3.3)). 
In sections 4, 5 and 6, we prove 
these three propositions. 
\par
\vskip 4pt
\noindent
{\bf Acknowledgement.}
This note has been grown up from my joint project [IOZ] (in progress) 
and much inspired by previous works [Ko1, 2] and [OZ]. I would like to express 
my thanks to  Professors A. Ivanov and De-Qi Zhang for their several valuable 
discussions. I am supported in part by Grant-in-Aid for Scientific Research, 
Ministry of Education, Science and Culture, Japan.

\section{Some basic properties of K3 groups after Nikulin and Mukai}

Let $X$ be a K3 surface and $G$ be a K3 group acting on $X$ by $\rho : 
G \times X \longrightarrow X$. Then $G$ has a $1$-dimensional 
representation on $H^{0}(X, \Omega_{X}^{2}) 
= \mathbf C \omega_{X}$ by $g^{*}\omega_{X} = \alpha(g)\omega_{X}$, and 
we have 
the exact sequence, called the basic sequence: 
$$1 \longrightarrow G_{N} := \text{Ker}\, \alpha \longrightarrow G \mapright{\alpha} 
\mu_{I} \longrightarrow 1\,\,  .$$ 
We call $G_{N}$ the symplectic part and $\mu_{I}$ (resp. $I$) the 
transcendental 
part (resp. the transcendental value) of the action 
$\rho : G \times X \longrightarrow X$.

By the basic sequence, the study of K3 group is divided into 
three parts: study of symplectic K3 groups 
$G_{N}$, study of transcendental values $I$, and study 
of possible extensions of 
symplectic parts 
by transcenental parts from the right. 

\begin{example} \label{example:fermatI}
The group $\tilde{F}_{384} = (\mu_{4}^{4} : S_{4})/\mu_{4}$ fits into the following exact sequence: 
$$1 \longrightarrow \mu_{4}^{4}/\mu_{4} \longrightarrow \tilde{F}_{384} \mapright{p} S_{4} 
\longrightarrow 1\,\, .$$ 
Then the group $\langle (1324), (34) \rangle \simeq D_{8}$ is a (one of three) Sylow $2$-subgroup of $S_{4}$ 
and $p^{-1}(\langle (1324),  (34) \rangle)$ is a (one of three) $2$-Sylow subgroup of 
$\tilde{F}_{384}$. We fix $\tilde{F}_{128}$ as this subgroup. 
The basic sequences of the standard actions by 
$\tilde{F}_{384}$ and by $\tilde{F}_{128}$ on the Fermat K3 surface 
$X_{4}$ are as follows:
$$1 \longrightarrow F_{384} := (\tilde{F}_{384})_{N} \longrightarrow 
\tilde{F}_{384} \mapright{\alpha} \mu_{4} 
\longrightarrow 1\,\, ,$$ 
$$1 \longrightarrow F_{128} := (\tilde{F}_{128})_{N}  \longrightarrow \tilde{F}_{128} \mapright{\alpha} \mu_{4} 
\longrightarrow 1\,\, .$$ 
The orders of the symplectic parts $F_{384}$ and $F_{128}$ are 
$384 = 2^{7}\cdot 3$ and $128 = 2^{7}$ respectively. 
Moreover, both basic sequences split: 
$\tilde{F}_{384} = F_{384} : \mu_{4}$ and 
$\tilde{F}_{128} = F_{128} : \mu_{4}$. Here 
the splittings are given by $\alpha(\text{diag}\,(1, 1, 1, \zeta_{4})) 
= \zeta_{4}$. 
\end{example}  
The next theorem by Nikulin [Ni] is the first important result 
about the symplectic part: 

\begin{theorem} \label{thm:numberTh} [Ni]
Let $g \in G_{N}$. Then $\text{ord}\, g \le 8$. The fixed locus 
$X^{g}$ is a finite set (if $g \not= 1$) 
and the cardinality $\vert X^{g} \vert$ depends only on $\text{ord}\, g$ 
as in the following table: 
$$ 
\begin{tabular}{*{9}{|c}{|}} 
\hline 
$ord(g)$ & 1 & 2 & 3 & 4 & 5 & 6 & 7 & 8 \\ 
\hline 
$\vert X^{g} \vert$ & $X$ & 8 & 6 & 4 & 4 & 2 & 3 & 2 \\
\hline 
\end{tabular} 
$$
\end{theorem}

Let $\Omega := \{1, 2, \cdots , 24\}$ be the set of $24$ elements and 
$\mathcal P(\Omega)$ be the power set of $\Omega$, i.e. the set of all 
subsets of $\Omega$. As it is classically known 
(see for instance [CS, Chap. 10]), $\mathcal P(\Omega)$ has a very 
remarkable subset $St(5, 8, 24)$, 
called the Steiner system. $St(5, 8, 24)$ is defined to be a subset of 
$\mathcal P(\Omega)$ consisting 
of $8$-element subsets such that for each $5$-element subset $B$ of 
$\Omega$, there is exactly one 
$A \in St(5, 8, 24)$ such that $B \subset A$. Such subsets 
$St(5, 8, 24)$ of $\mathcal P(\Omega)$ 
are known to be unique up to 
$\text{Aut}\, \Omega = S_{24}$ and satisfy $\vert St(5, 8, 24) \vert = 759$. We fix one of such $St(5, 8, 24)$. The Mathieu group $M_{24}$ of degree $24$ is 
then defined to be the stabilizer group of 
$St(5, 8, 24)$:
$$M_{24} := \{\tau \in \text{Aut}(\Omega) = 
S_{24} \vert \tau(St(5, 8, 24)) = St(5, 8, 24) \}\,\, .$$
It is well known that $M_{24}$ is a simple (sporadic) group of order $2^{10}\cdot 3^{3} \cdot 5 \cdot 
7 \cdot 11 \cdot 23$ and $5$-transitively acts on $\Omega$ (eg. [CS, ibid]). 
The 
Mathieu group $M_{23}$ 
of degree $23$ is the stabilizer group of one point, say 
$24 \in \Omega$, i.e. 
$M_{23} := \{\tau \in M_{24} \vert \tau(24) = 24 \}$.   
$M_{23}$ is also a simple group and is of order 
$\vert M_{23} \vert = \vert M_{24} \vert/24 = 2^{7} \cdot 3^{2} \cdot 5 \cdot 
7 \cdot 11 \cdot 23$. By the definition, both $M_{24}$ and $M_{23}$ act 
naturally on $\Omega$.

Mukai [Mu] discovered the following very beautiful theorem: 

\begin{theorem} \label{thm:MathieuTh} [Mu, main theorem]
Let $K$ be a finite group. Then $K$ is a symplectic K3 group on some K3 
surface $X$ if and only if $K$ is isomorphic to a subgroup 
of $M_{23}$ 
having at least $5$-orbits on $\Omega$ (under the action induced by 
the action by 
$M_{23}$ on $\Omega$). Moreover, with respect to the inclusion as abstract 
groups, there are 
exactly $11$ maximal such groups. The groups $M_{20}$ and $F_{384}$ are 
the ones of 
two largest orders, which are $\vert M_{20} \vert = 960$ and 
$\vert F_{384} \vert = 384$. 
\end{theorem} 

Later, Xiao [Xi] 
and Kondo [Ko] gave alternative proofs respectively. 
In the course of proof, Xiao 
shows that there are exactly $80$ isomorphism classes of symplectic 
K3 groups (as abstract 
groups). In our proof of the main result (1.2), we shall also exploit 
an idea of Kondo's 
alternative proof (section 6).

We emphasize the following consequence:
\begin{corollary} \label{corollary:nilpgpI} $F_{128}$ is isomorphic to a 
Sylow $2$-subgroup of $M_{23}$.  
\end{corollary} 
\begin{proof}
By example ({\ref{example:fermatI}}) and theorem 
({\ref{thm:MathieuTh}}), we have $F_{128} < M_{23}$. 
Moreover, since $\vert F_{128} \vert = 2^{7}$ and $\vert M_{23} 
\vert = 2^{7}\cdot k$ ($(2, k) = 1$), 
the result follows from the Sylow theorem. 
\end{proof}

Next we recall basic properties of the transcendental part $\mu_{I}$ of a K3 group on $X$ from 
[Ni] and [MO]. By $\varphi(I)$, we denote the Euler function of $I$, i.e. 
$\varphi(I) = \vert \text{Gal}(\mathbf Q(\zeta_{I})/\mathbf Q) \vert$. 
Note that $\varphi(I)$ is even unless $I = 1, 2$.

As it is observed in [Ni], $X$ is projective if $I \geq 2$. Indeed, if 
$I \geq 2$, then the quotient surface $X/G$ 
is either a rational surface with quotient singularities or 
an Enriques surface with at worst 
Du Val singularities, both of which are projective. In the rest of section 2, 
we assume that $X$ is projective.

Let $\text{NS}(X)$ be the N\'eron-Severi lattice of $X$ and $T(X)$ the transcendental lattice, 
i.e. the orthogonal complement of $\text{NS}(X)$ in $H^{2}(X, \mathbf Z)$ with respect to the cup product: 
$$T(X) := \{x \in H^{2}(X, \mathbf Z) \vert (x, \text{NS}(X)) = 0\}\,\, .$$
Then, $\text{NS}(X) \oplus T(X)$ is a sublattice of finite index of $H^{2}(X, \mathbf Z)$ 
(by the projectivity of $X$). $T(X)$ is also the minimal primitive 
sublattice of $H^{2}(X, \mathbf Z)$ such that the scalar extension by $\mathbf C$ contains 
the class of $\omega_{X}$ (by the Lefschetz $(1,1)$-theorem). Since $b_{2}(X) = 22$, we have $2 \leq \text{rank}\, T(X) \leq 21$. 

\begin{theorem} \label{thm:trans} ([Ni], see also [MO] for (2), (4)) 
\begin{list}{}{
\setlength{\leftmargin}{10pt}
\setlength{\labelwidth}{6pt}
}
\item[(1)] $G_{N}$ acts on $T(X)$ as identity. 
\item[(2)] Set $G/G_{N} = \langle g \, \text{mod}\, G_{N} \rangle 
\simeq \mu_{I}$. Then, there is a natural isomorphism
$$T(X) \simeq \mathbf Z[\zeta_{I}]^{\oplus n}\,\, ,\,\, n = 
\frac{\text{rank}\, T(X)}{\varphi(I)}$$ 
as $\mathbf Z[\zeta_{I}]$-modules. Here, $\mathbf Z[\zeta_{I}]$-module structure on $T(X)$ is given by $f(\zeta_{I})x := f(g^{*})x$.
\item[(3)] $\varphi(I) \vert \text{rank}\, T(X)$. In particular, $\varphi(I) \leq 20$ and $I \leq 66$. Moreover, $I = 1, 2, 3, 4, 6$ if $\varphi(I) \leq 2$ 
and $I = 5, 8, 10, 12$ if $\varphi (I) = 4$. 
\item[(4)] $I \not= 60$. Conversely, each $I$ such that $\varphi(I) \leq 20$ 
and $I \not= 60$ 
is realized as a transcendental value of some K3 group. There are exactly 
$40$ such $I$. (For the explicit list, see [MO].)
\end{list}
\end{theorem} 

As we reviewed above, both symplectic part and transcendental part are now well understood. 
However, the K3 groups, i.e. all the geometrically possible extensions of 
$80$ symplectic parts by $40$ transcendental parts, are not yet classified 
completely. 
This problem is now in progress in the work [IOZ] by A. Ivanov, De-Qi Zhang and the author. 

We close section 2 by recalling the following group theoretical nature of 
$F_{128}$ from [Xi]: 
\begin{proposition} \label{prop:nilpII}
\begin{list}{}{
\setlength{\leftmargin}{10pt}
\setlength{\labelwidth}{6pt}
}
\item[(1)] The order structure of $F_{128}$ is as follows:  
$$
\begin{tabular}{*{5}{|c}{|}} 
\hline 
\text{order} & 1 & 2 & 4 & 8\\ 
\hline 
\text{cardinality}  & 1 & 35 & 76 & 16\\
\hline 
\end{tabular} 
$$
\item[(2)] The commutator subgroup $[F_{128}, F_{128}]$ of $F_{128}$ is 
isomorphic to 
$C_{2} \times D_{8}$, where $C_{n}$ is a cyclic group of order $n$ and $D_{2n}$ is a 
dihedral group of order $2n$. 
\item[(3)] $F_{128}$ has a subgroup isomorphic to the binary dihedral group of order $16$: 
$$Q_{16} := \langle a, b \vert a^{8} = 1, a^{4} = b^{2}, b^{-1}ab = a^{-1} 
\rangle\,\, .$$ 
\end{list}
\end{proposition}
\begin{proof}
One can read off all these informations from Xiao's table [Xi]. Of course, one can directly 
check each 
statement from explicit descriptions of $\tilde{F}_{128}$ and $F_{128}$ in 
example ({\ref{example:fermatI}}). For instance, 
we have the following descriptions in $\text{PGL}(4, \mathbf C)$:
$$F_{128} = \{[(\zeta_{4}^{a_{i}}\delta_{i\sigma (j)})] \vert 
\sigma \in \langle (1324), (34) 
\rangle (\simeq D_{8}), sgn(\sigma)\cdot \Pi_{i=1}^{4}\zeta_{4}^{a_{i}} 
= 1\}$$ 
$$[F_{128}, F_{128}] = \langle A, B, C \rangle \simeq C_{2} \times D_{8}\,\, ;\,\, 
F_{128} > \langle P, Q \rangle \simeq Q_{16}\,\, ,$$ 
where 
$$A = \begin{pmatrix}
1&0&0&0\cr 
0&1&0&0\cr
0&0&-1&0\cr
0&0&0&-1
\end{pmatrix}\,\, ,\,\, 
B = \begin{pmatrix}
1&0&0&0\cr 
0&\zeta_{4}&0&0\cr
0&0&1&0\cr
0&0&0&\zeta_{4}^{-1}
\end{pmatrix}\,\,  ,\,\, 
C =  \begin{pmatrix}
0&1&0&0\cr 
1&0&0&0\cr
0&0&0&1\cr
0&0&1&0
\end{pmatrix}\,\, ,$$ 
$$
P = \begin{pmatrix}
0&1&0&0\cr 
\zeta_{4}&0&0&0\cr
0&0&1&0\cr
0&0&0&\zeta_{4}
\end{pmatrix}\,\,  ,\,\, 
Q = \begin{pmatrix}
\zeta_{4}&0&0&0\cr 
0&\zeta_{4}&0&0\cr
0&0&0&1\cr
0&0&1&0
\end{pmatrix}\,\, .$$
All the statements also follow from these descriptions. 
\end{proof}

\section{Reduction of the main theorem to three propositions}

In this section, we reduce the main theorem (1.2) to the following three 
propositions: 

\begin{proposition} \label{prop:propIIb} 
Let $X$ be a projective K3 surface. Assume that $Q_{16}$ is a 
symplectic K3 group on 
$X$. Then 
\begin{list}{}{
\setlength{\leftmargin}{10pt}
\setlength{\labelwidth}{6pt}
}
\item[(1)] 
$\text{NS}(X)^{Q_{16}} = \mathbf Z H$, where $H$ is an ample class on $X$. 
\item[(2)] If, in addition, $(H^{2}) = 4$, then the polarized K3 surface 
$(X, H)$ is unique up to 
isomorphism. In particular, $(X, H) \simeq 
(X_{4}, H_{4})$, where $X_{4}$ is the Fermat quartic K3 surface and 
$H_{4} := \iota^{*} \mathcal 
O_{\mathbf P^{3}}(1)$ under the natural inclusion $\iota: X_{4} 
\subset \mathbf P^{3}$. 
\end{list}
\end{proposition}

\begin{proposition} \label{prop:propIIa}
\begin{list}{}{
\setlength{\leftmargin}{10pt}
\setlength{\labelwidth}{6pt}
}
\item[(1)] Let $G$ be a K3 group on $X$ such that $G_{N} \simeq F_{384}$. 
Then the transcendental value $I$ of $G$ is either $1$, $2$, or $4$.   
\item[(2)] Let $G$ be a solvable K3 group on $X$. Then $\vert G \vert \leq 2^{9}\cdot 3$. 
Moreover, if $\vert G \vert = 2^{9} \cdot 3$, then the symplectic part $G_{N}$ 
is necessarily isomorphic to $F_{384}$ and the transcendental part is isomorphic to 
$\mu_{4}$. 
\item[(3)] Let $G$ be a K3 group on $X$ such that $G_{N} \simeq F_{128}$. 
Then the transcendental value $I$ of $G$ is either $1$, $2$, or $4$. In 
particular, $G$ is a $2$-group and nilpotent.   
\item[(4)] Let $G$ be a nilpotent K3 group on $X$. Then $\vert G \vert \leq 2^{9}$. 
Moreover, if $\vert G \vert = 2^{9}$, then the symplectic part $G_{N}$ 
is necessarily isomorphic to $F_{128}$ and the transcendental part is isomorphic to 
$\mu_{4}$. 
\end{list}
\end{proposition}
 
\begin{proposition} \label{prop:propIIc}
Let $X$ be a K3 surface. Assume that $X$ admits a K3 group $G$ of order 
$2^{9}$. 
Then $X$ is projective and $\text{NS}(X)^{G} = \mathbf Z H$, where 
$H$ is an ample class such that $(H^{2}) = 4$. 
\end{proposition}

We shall prove these three propositions in sections 4, 5, 6 respectively. 
In the rest of this section, 
we show that these propositions imply the main result (1.2). 
\vskip0.5cm
\noindent
{\it Proof of "(3.1), (3.2) and (3.3) imply (1.2)"}
\vskip0.3cm

Let $Y$ be a K3 surface admitting a K3 group $F$ such that $\vert F 
\vert = 2^{9} \cdot 3$. 
Let $G$ be a Sylow $2$-subgroup of $F$. Then $\vert G \vert = 2^{9}$ and 
$G$ is a nilpotent group. 
(Here we recall that any $p$-group is nilpotent.) Then, by (3.2)(4), 
$G_{N} \simeq F_{128}$ and $I = 4$. In particular, $Y$ is projective by 
$I \geq 2$. Recall that  
$Q_{16}$ is a subgroup of $F_{128}$ by (2.6)(3). Then, we have embeddings: 
$Q_{16} < F_{128} < G < F$.  
Thus 
$$\text{NS}(Y)^{F} \subset \text{NS}(Y)^{G} \subset 
\text{NS}(Y)^{F_{128}} \subset 
\text{NS}(Y)^{Q_{16}} = \mathbf Z H\,\, .$$ 
Here we use (3.1)(1) for the last equality. Since 
$\text{NS}(Y)^{F}$ contains 
an ample invariant class, say $\sum_{g \in F} g^{*}h$, $h$ being ample 
on $Y$, we have then:
$$ \text{NS}(Y)^{F} = \text{NS}(Y)^{G} = \text{NS}(Y)^{F_{128}} = 
\text{NS}(Y)^{Q_{16}} = 
\mathbf Z H\,\, .$$
Hence $F < \text{Aut}(X, H)$. Moreover, $(H^{2}) = 4$ by (3.3). Hence, 
$\varphi : (X, H) \simeq (X_{4}, H_{4})$ by (3.1)(2). Then, under the 
isomorphism $F \simeq 
\varphi^{-1} \circ F \circ \varphi$, we have $((X, H), F) \simeq 
((X_{4}, H_{4}), F)$. So, we may identify $((X, H), F) = 
((X_{4}, H_{4}), F)$. Under this identification, we have 
$$F < \text{Aut}(X_{4}, H_{4}) > \tilde{F}_{384} > F_{384}\,\,  .$$
Note that $\text{Aut}(X_{4}, H_{4})$ is a finite group. This is because 
$\text{Aut}(X_{4}, H_{4})$ is a discrete algebraic subgroup of 
$PGL(\mathbf P^{3})$, whence, finite.  Thus, $[\text{Aut}(X_{4}, H_{4}) 
: F_{384}] \leq 4$ by (3.2)(1) and (2.3). Hence $\vert \text{Aut}(X_{4}, H_{4}) \vert 
\leq 2^{7}\cdot 3 \times 4 = 2^{9}\cdot 3$. Since 
$\vert \tilde{F}_{384} \vert = \vert F \vert = 2^{9} \cdot 3$, we then obtain  
$F = \text{Aut}(X_{4}, H_{4}) = \tilde{F}_{384}$.  
This implies the assertion (1) of the main theorem.

Next, we shall show the assertion (2) of the main theorem. Let $X$ be a K3 surface admitting 
a K3 group $G$ such that $\vert G \vert = 2^{9}$. Then, by repeating the same argument as 
above, we can identify 
$((X, H), G) = ((X_{4}, H_{4}), G)$. Since 
$\text{Aut}(X_{4}, H_{4}) = \tilde{F}_{384}$, our 
$G$ is a subgroup of $\tilde{F}_{384}$. Since $\vert \tilde{F}_{384} \vert = 2^{9} \cdot 3$ and 
$\vert G \vert = 2^{9}$, it follows that $G$ is one of three Sylow $2$-subgroups $\tilde{F}_{128}$ 
of $\tilde{F}_{384}$, which are conjugate one another in 
$\tilde{F}_{384}$. This implies the result.

\begin{remark}\label{remark:shioda} As a byproduct, we have 
obtained that
$\text{Aut}(X_{4}, H_{4}) = \tilde{F}_{384}$. One can also derive this equality by a more direct calculation along the same line as in [Sh]. We also 
notice that $\rho(X_{4}) = 20$ and, by [SI], the full automorphism group 
$\text{Aut}(X_{4})$ is an infinite group. 
\end{remark}

\section{Polarized K3 surface of degree $4$ with a symplectic $Q_{16}$-action}

In this section, we shall prove proposition (3.1). 
\begin{definition}\label{definition:binary} The binary dihedral group $Q_{4m}$ of order $4m$ 
is defined by
$$Q_{4m} := \langle a, b \vert a^{2m} = 1, a^{m} = b^{2}, b^{-1}ab = a^{-1} \rangle\,\, .$$ 
$Q_{4m}$ is realized as a linear subgroup of $\text{GL}(2, \mathbf C)$ as 
$$\langle a :=  \begin{pmatrix}
\zeta_{2m}&0\cr 
0&\zeta_{2m}^{-1}
\end{pmatrix}\,\,  ,\,\, 
b := \begin{pmatrix}
0&\zeta_{4}\cr 
\zeta_{4}&0
\end{pmatrix}\,\, \rangle  = \{a^{n}, a^{n}b\, \vert\, 0 \leq n \leq 2m -1 \}.$$ 
\end{definition}

The next two lemmas (4.2) and (4.4)($=$(3.1)(1)) explain the reason why we 
pay a special attension to the particular group $Q_{16}$.

\begin{lemma} \label{lemma:linear} Any projective representation of $Q_{4m}$ is induced 
by a linear representation, i.e. for any group homomorphism
$\rho : Q_{4m} \longrightarrow PGL(n, \mathbf C) := GL(n, \mathbf C)/\mathbf C^{\times}$,  
there is a group homomorphism $\tilde{\rho} : Q_{4m} \longrightarrow GL(n, \mathbf C)$ 
such that $\rho = p \circ \tilde{\rho}$, where $p : GL(n, \mathbf C) \longrightarrow PGL(n, \mathbf C)$ 
is the quotient map. 
\end{lemma} 

\begin{proof} We shall write 
$[X] = X \text{mod}\, \mathbf C^{\times}$ for $X \in \text{GL}(n, \mathbf C)$. 
First we remark that $Q_{4m} = \langle a, b \vert a^{m} = b^{2}, b^{-1}ab = a^{-1} \rangle$.  
This is because $a^{m} = b^{2}$ and $b^{-1}ab = a^{-1}$ imply that $a^{-m} = b^{-1}a^{m}b 
= b^{-1}b^{2}b = b^{2} = a^{m}$, whence $a^{2m} = 1$.  Let $\rho(a) = [A]$ and $\rho(b) = [B]$. 
Then, $[A^{m}] = [B^{2}]$ and $[B^{-1}AB] = [A^{-1}]$. I.e. $A^{m} = \alpha B^{2}$ and 
$B^{-1}AB = \beta A^{-1}$ in $\text{GL}(n, \mathbf C)$ for some $\alpha, \beta \in 
\mathbf C^{\times}$. By replacing the representative $A$ by $A/\sqrt{\beta}$, 
one has 
$B^{-1}AB = A^{-1}$. Next, by replacing the representative $B$ by 
$\sqrt{\alpha}B$, 
one obtains $B^{-1}AB = A^{-1}$ and $A^{m} = B^{2}$. Therefore we have a group 
homomorphism $\tilde{\rho} : Q_{4m} \longrightarrow \text{GL}(n, \mathbf C)$ defined by 
$\tilde{\rho}(a) = A$ and $\tilde{\rho}(b) = B$. This $\tilde{\rho}$ satisfies 
$\rho = p \circ \tilde{\rho}$. 
\end{proof}
\begin{remark} \label{remark:nonlinear} Consider 
 the dihedral group $D_{8} := \langle a, b \vert 
 a^{4} = b^{2} = 1, b^{-1}ab = a^{-1} \rangle$. Then the map 
$$\rho(a) = \begin{pmatrix}
\zeta_{8}&0\cr 
0&\zeta_{8}^{-1}
\end{pmatrix}\,\,  ,\,\, 
\rho(b) = \begin{pmatrix}
0&\zeta_{4}\cr 
\zeta_{4}&0
\end{pmatrix}$$ 
defines a projective representation 
$\rho : D_{8} \longrightarrow PGL(2, \mathbf C)$. However, this is 
not induced by any linear 
representation $D_{8} \longrightarrow GL(2, \mathbf C)$.  
\end{remark}

\begin{lemma} \label{lemma:polarization} Let $X$ be a projective K3 surface 
admitting 
a symplectic K3 group $Q_{16}$. Then $\text{rank}\, T(X) = 2$ 
and $\text{NS}(X)^{Q_{16}}  = \mathbf Z H$ 
for some ample 
class $H$. 
\end{lemma}

\begin{proof}  The following very useful lemma is essentially 
due to Mukai [Mu]. 

\begin{proposition} \label{lemma:polarization} Let $X$ be a projective K3 surface admitting 
symplectic K3 group $G$. Then 
\begin{list}{}{
\setlength{\leftmargin}{10pt}
\setlength{\labelwidth}{6pt}
}
\item[(1)] 
$$\text{rank}\, H^{2}(X, \mathbf Z)^{G}  = \frac{1}{\vert G \vert}(24 + \sum_{n =2}^{8} 
m(n)f(n)) - 2\,\, ,$$ 
where $m(n)$ is the number of elements of order $n$ in $G$ and $f(n)$ is the number 
of the fixed 
points in (2.2). 
\item[(2)] $\text{rank}\, H^{2}(X, \mathbf Z)^{G}  \geq 3$. Moreover, if 
$\text{rank}\, H^{2}(X, \mathbf Z)^{G}  = 3$, 
then $\text{rank}\, T(X) = 2$ and $\text{NS}(X)^{G} = \mathbf Z H$ for some 
ample class $H$.  
\end{list}
\end{proposition}

\begin{proof} Consider the action by $G$ on the total cohomology group 
$H^{*}(X, \mathbf Z) = 
H^{0}(X, \mathbf Z) \oplus H^{2}(X, \mathbf Z) \oplus 
H^{4}(X, \mathbf Z)$. Then, by the representation theory, one has
$$\text{rank}\, H^{*}(X, \mathbf Z)^{G} = 
\frac{1}{\vert G \vert}\sum_{g \in G} 
\text{tr}(g^{*} \vert H^{*}(X, \mathbf Z))\,\, .$$ 
By the Lefschetz fixed point formula, each summand satisfies  
$$\text{tr}(g^{*} \vert H^{*}(X, \mathbf Z)) = e(X^{g})\,\, .$$ 
Here $e(*)$ is the topological Euler number of $*$. 
Now, by (2.2) and $e(X) = 24$, we have 
$$\text{rank}\, H^{*}(X, \mathbf Z)^{G} =  
\frac{1}{\vert G \vert}(24 + \sum_{n =2}^{8} m(n)f(n))\,\, 
.$$ 
Since $G$ acts trivially on $H^{0}(X, \mathbf Z) \simeq H^{4}(X, \mathbf Z) \simeq \mathbf Z$, 
we have the assertion (1).  Since $X$ is projective, $X$ has an invariant 
ample class. On the other hand, since $G$ is symplectic, 
$T(X)$ is a subset of $H^{2}(X, \mathbf Z)^{G}$ and is of rank at least $2$. 
Thus, $\text{rank}\, H^{2}(X, \mathbf Z)^{G} \geq 3$, and if the equality 
holds, 
then we have 
$\text{rank}\, T(X) = 2$ and $NS(X)^{G} = \mathbf Z H$ for some ample class 
$H$. This completes the proof. 
\end{proof} 
We also remark the following: 
\begin{lemma} \label{lemma:order} The order 
structure of $Q_{16}$ is as follows:  
$$
\begin{tabular}{*{5}{|c}{|}} 
\hline 
\text{order} & 1 & 2 & 4 & 8\\ 
\hline 
\text{cardinality}  & 1 & 1 & 10 & 4\\
\hline 
\end{tabular} 
$$
\end{lemma}
\begin{proof} This directly follows from the description of $Q_{16}$. 
\end{proof} 

Let us return back to (4.4). By (4.5) and (4.6), we calculate that 
$$\text{rank}\, H^{2}(X, \mathbf Z)^{Q_{16}}  = 
\frac{1}{16}(24 + 8\cdot 1 + 4 \cdot 10 + 2 \cdot 4 ) - 2 = 3\,\, .$$ 
This completes the proof of (4.4). 
\end{proof} 

Recall that the standard action by $F_{128}$ on $X_{4}$ is a symplectic 
action on the polarized 
K3 surface $(X_{4}, H_{4})$ and $Q_{16} < F_{128}$. Now, the next proposition 
completes 
the proof of proposition (3.1): 

\begin{proposition} \label{proposition:qaurtic} Polarized K3 surfaces $(X, H)$ of degree $4$ which 
admit a symplectic K3 group $Q_{16}$ (which keeps $H$ invarinat) are 
unique up to isomorphism as 
polarized K3 surfaces.  
\end{proposition}
\begin{proof} Since $\text{rank}\, \text{NS}(X)^{Q_{16}} = 1$ and 
$H \in \text{NS}(X)$ 
is primitive by $(H^{2}) = 4$ and by even-ness of the intersection numbers, 
we have 
$\text{NS}(X)^{Q_{16}} = \mathbf Z H$. Then $\vert H \vert$ has no fixed 
components. Indeed, the fixed part of $\vert H \vert$ must be also 
$Q_{16}$-invariant, while $\text{NS}(X)^{Q_{16}} 
= \mathbf Z H$. 
Therefore, the ample linear system $\vert H \vert$ is free by [SD]. 
Note that $\text{dim} \vert H \vert = 3$ by
the Riemann-Roch formula and by $(H^{2}) = 4$. Then $\vert H \vert$
defines a morphism 
$$\Phi := \Phi_{\vert H \vert} : X \longrightarrow 
\mathbf P^{3} = \vert H \vert^{*}\,\, ;\,\, x \mapsto \{D \in \vert H \vert\, 
\vert\, D \ni x \, \}\,\, .$$
This $\Phi$ is either

\begin{list}{}{
\setlength{\leftmargin}{10pt}
\setlength{\labelwidth}{6pt}
}
\item[(I)] an embedding onto a (smooth) quartic surface $W =(4)$, or 
\item[(II)] a finite double cover of an irreducible, reduced quadratic surface 
$W =(2)$. 
\end{list}

Since $H$ is $Q_{16}$-invariant in $\text{Pic}(X) \simeq \text{NS}(X)$, 
the divisor $g^{*}D$ is linearly equivalent to 
$D$ whenever $D \in \vert H \vert$ and $g \in Q_{16}$. Thus, 
the group $Q_{16}$ induces a $\Phi$-equivariant, projective linear action 
on the image $W$. 
By (4.2), this action is also induced by a 
linear co-action by $Q_{16}$ on
$H^{0}(X, \mathcal O_{X}(H)) = \oplus_{i=1}^{4} \mathbf Cx_{i}$.  

In order to complete the proof, it suffices to show the two assertions, 
that the case (II) can not happen (4.9) and 
that the image $W$ is uniquely determined up to projective transformations of 
$\mathbf P^{3}$ in the case (I) (4.11). 

In both assertions, we need the following classification of the complex irreducible 
linear representations of $Q_{16} = \langle a, b \vert a^{4} = b^{2}, b^{-1}ab = a^{-1} \rangle$: 
\begin{lemma} \label{lemma:irredrep} A complex irreducible 
linear representation of 
$Q_{16}$ is isomorphic to one of 
the following $7$ representations: 
$$ \rho_{1, 1} : a \mapsto 1, b \mapsto 1\,\, ,\,\, \rho_{1, 2} : a \mapsto 1, b \mapsto -1\,\, ,$$ 
$$ \rho_{1, 3} : a \mapsto -1, b \mapsto 1\,\, ,\,\, \rho_{1, 4} : a \mapsto -1, b \mapsto -1\,\, ,$$ 
$$\rho_{2, 1} : a \mapsto \begin{pmatrix}
\zeta_{8}&0\cr 
0&\zeta_{8}^{-1}
\end{pmatrix}\,\,  ,\,\, 
b \mapsto \begin{pmatrix}
0&\zeta_{4}\cr 
\zeta_{4}&0
\end{pmatrix}\,\,  ,\,\, \rho_{2,2} : a \mapsto \begin{pmatrix}
\zeta_{8}^{3}&0\cr 
0&\zeta_{8}^{-3}
\end{pmatrix}\,\,  ,\,\, 
b \mapsto \begin{pmatrix}
0&\zeta_{4}\cr 
\zeta_{4}&0
\end{pmatrix}\,\,  ,$$
$$\rho_{2, 3} : a \mapsto \begin{pmatrix}
\zeta_{4}&0\cr 
0&\zeta_{4}^{-1}
\end{pmatrix}\,\,  ,\,\, 
b \mapsto \begin{pmatrix}
0&1\cr 
1&0
\end{pmatrix}\,\,  .$$ 
\end{lemma} 

\begin{proof} These $7$ representations are clearly irreducible and well-defined. Moreover, 
any two are inequivalent as linear representations 
(by looking at the trace of the matrices). Since 
$16 = 1^{2} \cdot 4 + 2^{2} \cdot 3$, these are all. \end{proof} 

\begin{lemma} \label{lemma:quadric} The case (II) can not happen. \end{lemma} 
 
\begin{proof} In what follows, assuming to the contrary that the case (II) happens, i.e. 
the image $W$ is a quadratic surface, we shall 
derive a contradiction. 

\begin{claim} \label{claim:smooth} 
\begin{list}{}{
\setlength{\leftmargin}{10pt}
\setlength{\labelwidth}{6pt}
}
\item[(1)] $W$ is non-singular.
\item[(2)] The induced action by $Q_{16}$ on $W$ is faithful. 
\end{list}
\end{claim} 

\begin{proof} Note that a quadratic surface is normal if it is irreducible 
and reduced. Since $\Phi$ is a finite double covering, it is also a Galois covering. 
Let $\tau$ be  the covering involution. Then $W = X/\tau$. Since $W$ is a rational surface, 
$\tau^{*}\omega_{X} = -\omega_{X}$. Thus, if $P \in X^{\tau}$, then there is a local coordinate 
$(x_{P}, y_{P})$ at $P$ such that $\tau^{*}(x_{P}, y_{P}) = (x_{P}, -y_{P})$. 
Hence, $W$ 
is non-singular. The kernel of the natural map 
$\text{Aut}(X, H) \longrightarrow \text{Aut}(W, \mathcal O_{W}(1))$ is a subgroup of $\langle \tau \rangle$. 
Since $\tau^{*}\omega_{X} = -\omega_{X}$, we have $Q_{16} \cap \langle \tau \rangle = \{1\}$. 
This means that the induced action by $Q_{16}$ on $W$ is faithful.  \end{proof} 

Let us return back to the proof of (4.9). Let us consider the irreducible decomposition 
of the co-action by $Q_{16}$ on 
$H^{0}(X, \mathcal O_{X}(H)) \simeq \mathbf C^{4}$. 
Note that the representations $\rho_{2,1}$ and $\rho_{2,2}$ are 
transformed by the outer automorphism $a \mapsto a^{3}$ and 
$b \mapsto b$ of $Q_{16}$. Recall also that the action must be faithful 
by (4.10)(2). Thus, we may assume without loss of generality that 
$\rho_{2,1}$ appears in the decomposition. Under this 
assumption, 
there are $4$ possible decompositions: (i) $\rho_{2,1} \oplus \rho_{2,1}$, 
(ii) $\rho_{2,1} \oplus \rho_{2,2}$, (iii) $\rho_{2,1} \oplus \rho_{2,3}$, 
(iv) 
$\rho_{2,1} \oplus (\text{two}\, 1-\text{dimensional 
irreducible representations})$. 
In cases (i) and (ii), $a^{4} = id$ in $\text{PGL}(4, \mathbf C)$, 
a contradiction to (4.10)(2).

Consider the case (iii). Then the action by $Q_{16}$ on 
$H^{0}(X, \mathcal O_{X}(H))$ is 
given by 
$$a =  \begin{pmatrix}
\zeta_{8}&0&0&0\cr 
0&\zeta_{8}^{-1}&0&0\cr
0&0&\zeta_{4}&0\cr
0&0&0&\zeta_{4}^{-1}
\end{pmatrix}\,\,  ,\,\, 
b = \begin{pmatrix}
0&\zeta_{4}&0&0\cr 
\zeta_{4}&0&0&0\cr
0&0&0&1\cr
0&0&1&0
\end{pmatrix}\,\, ,$$ 
under suitable basis  $\langle x_{i} \rangle_{i=1}^{4}$ of $H^{0}(X, \mathcal O_{X}(H))$. Let us 
consider the defining equation $F_{2}(x_{1}, x_{2}, x_{3}, x_{4}) \in 
\text{Sym}^{2}H^{0}(X, \mathcal O_{X}(H))$ of $W$. Then, $F_{2}$ is both 
$a$-semi-invariant and 
$b$-semi-invariant, i.e. $a(F_{2}) = \sigma(a)F_{2}$ and $b(F_{2}) = \sigma(b)F_{2}$. This $\sigma$ 
defines a $1$-dimensional representation of $Q_{16}$. Thus $a(F_{2}) = \pm F_{2}$ and 
$b(F_{2}) = \pm F_{2}$. If $a(F_{2}) = F_{2}$, then $F_{2} = \alpha x_{1}x_{2} + \beta x_{3}x_{4}$ 
by the explicit matrix form of $a$. Since $b(F_{2}) = \pm F_{2}$, we have  
then $F_{2} = \alpha x_{1}x_{2}$ 
or $\beta x_{3}x_{4}$. However, this contradicts the smoothness of $W$. 
If $a(F_{2}) = -F_{2}$, then $F_{2}$ is of the form $F_{2}(x_{3}, x_{4})$ and again contradicts 
the smoothness of $W$.  Thus, the case (iii) can not happen, either. 

Finally consider the case (iv). In this case, the action by $Q_{16}$ on $H^{0}(X, \mathcal O_{X}(H))$ is 
given by 
$$a =  \begin{pmatrix}
\zeta_{8}&0&0&0\cr 
0&\zeta_{8}^{-1}&0&0\cr
0&0&\pm 1&0\cr
0&0&0&\pm 1
\end{pmatrix}\,\,  ,\,\, 
b = \begin{pmatrix}
0&\zeta_{4}&0&0\cr 
\zeta_{4}&0&0&0\cr
0&0&\pm 1&0\cr
0&0&0&\pm 1
\end{pmatrix}\,\, ,$$ 
under suitable basis $\langle x_{i} \rangle_{i=1}^{4}$ of $H^{0}(X, \mathcal O_{X}(H))$. 
Let us consider the defining equation $F_{2}$ of $W$. Then as before 
$a(F_{2}) = \pm F_{2}$ and $b(F_{2}) = \pm F_{2}$. If $a(F_{2}) = -F_{2}$, then 
$F_{2} = F_{2}(x_{3}, x_{4})$ and $W$ is singular, a contradiction. Consider the case 
where $a(F_{2}) = F_{2}$. By the explicit form of $a$, we have 
$F_{2} = \alpha x_{1}x_{2} + f_{2}(x_{3}, x_{4})$. Since $W$ is non-singular, 
we have $\alpha \not= 0$ 
and $f_{2} \not= 0$. Since $b(x_{1}x_{2}) = -x_{1}x_{2}$, we have 
$b(f_{2}) = -f_{2}$. Thus, again by the explicit form of $b$, it follows that $F_{2} = \alpha x_{1}x_{2} 
+ \beta x_{3}x_{4}$ for some non-zero constants $\alpha$, $\beta$. After 
replacing $x_{i}$ by their multiples 
and the order of $x_{3}$ and $x_{4}$ if necessary, we finally 
normalize the equation of $W$ as $F_{2} = x_{1}x_{2} + x_{3}x_{4}$
and we have:
$$a =  \begin{pmatrix}
\zeta_{8}&0&0&0\cr 
0&\zeta_{8}^{-1}&0&0\cr
0&0&1&0\cr
0&0&0&1
\end{pmatrix}\,\, \text{or}\,\,   
\begin{pmatrix}
\zeta_{8}&0&0&0\cr 
0&\zeta_{8}^{-1}&0&0\cr
0&0&-1&0\cr
0&0&0&-1
\end{pmatrix}\,\,  ,\,\, \text{and}\,\, 
b = \begin{pmatrix}
0&\zeta_{4}&0&0\cr 
\zeta_{4}&0&0&0\cr
0&0&1&0\cr
0&0&0&-1
\end{pmatrix}\,\, .$$
Then, it follows that 
$W^{a} = W^{a^{2}} = W^{a^{4}} = \{P_{i}\}_{i=1}^{4} =: S$,  
where $P_{1} = [1 : 0: 0: 0]$, $P_{2} = [0:1:0:0]$, $P_{3} = [0:0:1:0]$ and $P_{4} = [0:0:0:1]$. 
Since the actions by $Q_{16}$ on $X$ and on $W$ are $\Phi$-equivariant and 
since $\Phi$ is  a finite 
morphism of degree $2$, it follows that $a^{2}$ and $a^{4}$ act on $T := \Phi^{-1}(S)$ 
as identity. Thus $X^{a^{2}} = X^{a^{4}} = T$. On the other hand, 
$\vert X^{a^{2}} \vert = 4$ 
and $\vert X^{a^{4}} \vert = 8$ by (2.2), a contradiction. This completes 
the proof of (4.9). 
\end{proof} 

\begin{lemma} \label{lemma:quartic} Assume that the case (I) happens, i.e. that 
$\Phi : X \simeq W = (4) \subset \mathbf P^{3}$. Then 
$W = (x_{1}^{4} + x_{2}^{4} + x_{3}^{3}x_{4} + x_{3}x_{4}^{3} = 0)$ 
in a suitably chosen homogeneous coordinates of $\mathbf P^{3}$.  \end{lemma} 

\begin{proof} Set $W = (F_{4}(x_{1}, x_{2}, x_{3}, x_{4}) = 0)$. 
We note that $\Phi$-equivariant action by $Q_{16}$ on $W$ is symplectic 
and faithful. 
As in the 
previous lemma (4.9), 
we consider the irreducible decomposition of the co-action by $Q_{16}$ on 
$H^{0}(X, \mathcal O_{X}(H))$. 
Again as before, we may assume that $\rho_{2,1}$ appears in the 
decomposition. Under this 
assumption, 
there are $4$ possible decompositions: (i) $\rho_{2,1} \oplus \rho_{2,1}$, 
(ii) $\rho_{2,1} \oplus \rho_{2,2}$, (iii) $\rho_{2,1} \oplus \rho_{2,3}$, (iv) 
$\rho_{2,1} \oplus (\text{two}\, 1-\text{dimensional irreducible 
representations})$. 
As before, the cases (i) and (ii) 
are ruled out by $a^{4} = id$ in $\text{PGL}(4, \mathbf C)$.

\begin{claim} \label{lemma:caseI} The case (iii) does not happen. \end{claim} 

\begin{proof} Assume that the case (iii) happens. Then the action by $Q_{16}$ on 
$H^{0}(X, \mathcal O_{X}(H))$ is 
given by 
$$a =  \begin{pmatrix}
\zeta_{8}&0&0&0\cr 
0&\zeta_{8}^{-1}&0&0\cr
0&0&\zeta_{4}&0\cr
0&0&0&\zeta_{4}^{-1}
\end{pmatrix}\,\,  ,\,\, 
b = \begin{pmatrix}
0&\zeta_{4}&0&0\cr 
\zeta_{4}&0&0&0\cr
0&0&0&1\cr
0&0&1&0
\end{pmatrix}\,\, ,$$ 
under suitable basis $\langle x_{i} \rangle_{i=1}^{4}$ of $H^{0}(X, \mathcal O_{X}(H))$. 
Since $\text{det}\, a = 1$ and $a^{*}\omega_{W} = \omega_{W}$, it follows that 
$F_{4}$ is $a$-invariant. By the explicit form of $a$, the 
equation $F_{4}$ must be then of the following form:
$$F_{4} = \alpha x_{1}^{2}x_{2}^{2} + \beta x_{1}x_{2}x_{3}x_{4} + f_{4}(x_{3}, x_{4})\,\, .$$ 
However the point $[1:0:0:0]$ is then a singular point of $W$, 
a contradiction. 
\end{proof}

In what follows, we shall consider the case (iv). Note that 
$$ \begin{pmatrix} 
\zeta_{8}&0&0&0\cr 
0&\zeta_{8}^{-1}&0&0\cr
0&0&-1&0\cr
0&0&0&-1
\end{pmatrix}\,\,  = \,\, 
\begin{pmatrix} 
\zeta_{8}^{5}&0&0&0\cr 
0&\zeta_{8}^{3}&0&0\cr
0&0&1&0\cr
0&0&0&1
\end{pmatrix}\,\,  = \,\, 
\begin{pmatrix}
\zeta_{8}&0&0&0\cr 
0&\zeta_{8}^{-1}&0&0\cr
0&0&1&0\cr
0&0&0&1
\end{pmatrix}^{5}$$ 
in $\text{PGL}(4, \mathbf C)$. So, replacing $a$ by $a^{5}$ by 
an outer automorphism 
$a \mapsto a^{5}$ and 
$b \mapsto b$ of $Q_{16}$ if necessary, we may assume that $a$ is either 
$$a_{1} :=  \begin{pmatrix} 
\zeta_{8}&0&0&0\cr 
0&\zeta_{8}^{-1}&0&0\cr
0&0&1&0\cr
0&0&0&1
\end{pmatrix}\,\, \text{or}\,\,  
a_{2} := \begin{pmatrix}
\zeta_{8}&0&0&0\cr 
0&\zeta_{8}^{-1}&0&0\cr
0&0&1&0\cr
0&0&0&-1
\end{pmatrix}\,\, .$$ 
In each case 
$$b = \begin{pmatrix} 
0&\zeta_{4}&0&0\cr 
\zeta_{4}&0&0&0\cr
0&0&\pm 1&0\cr
0&0&0& \pm 1
\end{pmatrix}\,\, .$$ 
\begin{claim} \label{lemma:caseII} $a \not= a_{1}$. \end{claim} 

\begin{proof} Assume that $a = a_{1}$. Then, as before, 
by $\text{det}\, a_{1} = 1$ and $a^{*}\omega_{W} 
= \omega_{W}$, it follows that $F_{4}$ is $a$-invariant. 
Thus $F_{4}$ must be of the following form: 
$$F_{4} = \alpha x_{1}^{2}x_{2}^{2} + \beta x_{1}x_{2}f_{2}(x_{3}, x_{4})  + 
f_{4}(x_{3}, x_{4})\,\, .$$ 
If $\beta = 0$, then $[1:0:0:0]$ is a singular point of $W$, a contradiction. 
If $\beta \not= 0$, we can normalize $\beta = 1$. If $\text{det}\, b = 1$, 
then $F_{4}$ is $b$-invariant by 
$b^{*}\omega_{X} = \omega_{X}$. Thus 
$b(f_{2}) = -f_{2}$ and then $f_{2} = 0$ by the explicit form of $b$. However, 
$[1:0:0:0]$ is then a singular point of $W$, a contradiction. 
If $\text{det}\, b = -1$, then $b(F_{4}) = -F_{4}$ by 
$b^{*}\omega_{X} = \omega_{X}$. Thus 
$\alpha  = 0$, $b(f_{2}) = f_{2}$ and $b(f_{4}) = -f_{4}$ by the explicit 
forms of $b$ and $F_{4}$. Then, $F_{4}$ is of the form 
$$F_{4} = x_{1}x_{2}(\beta_{1}x_{3}^{2} + \beta_{2}x_{4}^{2})  + \gamma_{1}x_{3}^{3}x_{4} 
+ \gamma_{2}x_{3}x_{4}^{3}\,\, .$$ 
However, the point $[0:1:0:0]$ is then a singular point of $W$, 
a contradiction. 
\end{proof} 

So, $a = a_{2}$. By $\text{det}\, a_{2} = -1$, we have $a(F_{4}) = -F_{4}$. By the 
explicit form of $a_{2}$, the equation $F_{4}$ 
is of the following form: 
$$F_{4} = \alpha x_{1}^{4} + \beta x_{2}^{4}  + \gamma x_{4}^{3}x_{3} + \delta x_{4}x_{3}^{3} 
+ \epsilon x_{1}x_{2}x_{3}x_{4}\,\, .$$

If $\alpha = 0$, then $\beta = 0$, because $F_{4}$ is $b$-semi-invariant. 
However, $[1:0:0:0]$ is then a singular point of $W$, a contradiction. 
Thus $\alpha \not= 0$. 
For the same reason, we have $\beta \not= 0$.  If $\gamma = 0$, then 
$[0:0:0:1]$ is a singular point of 
$W$. If $\delta = 0$, then $[0:0:1:0]$ is a singular point of $W$. Thus $\gamma \not= 0$ 
and $\delta \not= 0$.

If $\text{det}\, b = 1$, then $b(F_{4}) = F_{4}$. Thus 
$\alpha = \beta$ and $\epsilon = 0$. Then, applying a suitable 
linear transform like $x_{1} \mapsto cx_{1}$, $x_{2} \mapsto cx_{2}$, 
$x_{2} \mapsto dx_{2}$, 
$x_{3} \mapsto ex_{3}$, one can normalize the equation of $W$ as in (4.11).

If $\text{det}\, b = -1$, then $b(F_{4}) = -F_{4}$. Thus 
$\alpha = -\beta$ and $\epsilon = 0$. Then, one can again normalize 
the equation of $W$ as 
in (4.11). This completes the proof of (4.11). 
\end{proof}
We have now completed the proof of (4.7) and that of (3.1). 
\end{proof}

\section{Solvable K3 groups and nilpotent K3 groups}

In this section, we shall prove proposition (3.2). Throughout this section, 
we denote by $G$ a K3 group acting on $X$ and by 
$$1 \longrightarrow G_{N} \longrightarrow G\, \mapright{\alpha} \, 
\mu_{I} \longrightarrow 1\,\,$$ 
the basic sequence.

The next proposition is a special case of a more general fact in [IOZ] 
and is crucial for our proof: 

\begin{proposition} \label{proposition:split}
Assume that $I = 3$. Let $g$ be an element of $G$ such that $\alpha(g) = \zeta_{3}$.  
\begin{list}{}{
\setlength{\leftmargin}{10pt}
\setlength{\labelwidth}{6pt}
}
\item[(1)] Set $\text{ord}\, g = 3k$. Then $(k, 3) = 1$. In particular, the basic sequence splits if $I = 3$. 
\item[(2)] Assume that $\text{ord}\, g = 6$. Let $P \in X^{g}$. Then there is a local 
coordinate $(x, y)$ at $P$ such that either $g^{*}(x, y) = 
(\zeta_{6}^{-1}x, \zeta_{6}^{3}y)$ 
(type (1))  or  $g^{*}(x, y) = (\zeta_{6}^{-5}x, \zeta_{6}y)$ (type (5)). 
Let $m_{1}$, $m_{5}$ be the numbers of points of type (1) and of type 
(5). Then $(m_{1}, m_{5})$ 
is either $(2, 0)$, $(4,1)$ or $(6,2)$. 
\item[(3)] Assume that $\text{ord}\, g = 12$. Let $P \in X^{g}$. Then there is a local 
coordinate $(x, y)$ at $P$ such that either 
$g^{*}(x, y) = (\zeta_{12}^{-1}x, \zeta_{12}^{5}y)$ 
(type (1)), $g^{*}(x, y) = (\zeta_{12}^{-3}x, \zeta_{12}^{7}y)$ (type (3)), or 
$g^{*}(x, y) = (\zeta_{12}^{-9}x, \zeta_{12}y)$ (type (9)). 
Let $m_{1}$, $m_{3}$, $m_{9}$ be the numbers of points of types  (1), (3), 
(9) respectively. Then 
$(m_{1}, m_{3}, m_{9})$ is either $(1, 0, 0)$ or $(2, 1, 1)$. 
\end{list} 
\end{proposition}

\begin{proof} For the convenience to the readers, we shall give a 
proof for this special case. A more general treatment will be in [IOZ].

Let us show (1). If otherwise, $k = 3$ or $6$ by $g^{3} \in G_{N}$ and 
by (2.2). So, it suffices 
to show that $k \not= 3$.  Assume $ k = 3$. Then 
$\text{ord}\, g = 9$ and $X^{g} \subset 
X^{g^{3}}$. Since $X^{g^{3}}$ is a $6$-point set by (2.2), $X^{g}$ is 
also a finite set. 
Let $P \in X^{g}$. Then, since $\text{ord}\,g = 9$, $g^{*}\omega_{X} = 
\zeta_{3}\omega_{X}$, 
and $P \in X^{g}$ is isolated, there is a local 
coordinate $(x, y)$ at $P$ such that either 
$g^{*}(x, y) = (\zeta_{9}^{-1}x, \zeta_{9}^{4}y)$ 
(type (1)), $g^{*}(x, y) = (\zeta_{9}^{-2}x, \zeta_{9}^{5}y)$ (type (2)) or 
$g^{*}(x, y) = (\zeta_{9}^{-7}x, \zeta_{9}y)$ (type (7)) holds. 
Let $m_{1}$, $m_{2}$, $m_{7}$ be the numbers of fixed points of types 
(1), (2), (7). 
Then, by the holomorphic Lefschetz fixed point formula, one has: 
$$1 + \zeta_{3}^{-1} (= \sum_{i=0}^{2} (-1)^{i} 
\text{tr}(g^{*} \vert H^{i}(\mathcal O_{X}))$$ 
$$ = \frac{m_{1}}{(1- \zeta_{9}^{-1})(1 - \zeta_{9}^{4})} + 
\frac{m_{2}}{(1- \zeta_{9}^{-2})(1 - \zeta_{9}^{5})} + 
\frac{m_{7}}{(1- \zeta_{9}^{-7})(1 - \zeta_{9})}\,\, .$$ 
Note that the minimal polynomial of $\zeta_{9}$ over $\mathbf Q$ 
is $x^{6} + x^{3} + 1 = 0$. Now, a direct calculation shows that 
there is no solution 
$(m_{1}, m_{2}, m_{7})$ of the equation above even in $\mathbf Q$.

Let us give a proof of (2). For the same reason as in (1), one obtains 
a
list of possible local actions of $g$ at $P \in X^{g}$ as described in (2). 
Then, again by the holomorphic 
Lefschetz fixed point formula, one has: 
$$1 + \zeta_{3}^{-1} = \frac{m_{1}}{(1- \zeta_{6}^{-1})(1 - \zeta_{6}^{3})} + 
\frac{m_{5}}{(1- \zeta_{6}^{-5})(1 - \zeta_{6})}\,\, .$$ 
In addition, since $X^{g} \subset X^{g^{2}}$ and $X^{g^{2}}$ is a $8$-point set by (2.2), 
one has $m_{1} + m_{5} \leq 8$. Finding all 
the non-negative integer solutions $(m_{1}, m_{5})$ in this range, we obtain the result. Proof of (3) is similar. 
\end{proof}

The next lemma completes the assertions (1) and (3) of (3.2).  

\begin{lemma} \label{lemma:maximal} 
\begin{list}{}{
\setlength{\leftmargin}{10pt}
\setlength{\labelwidth}{6pt}
}
\item[(1)] If $G_{N} \simeq F_{128}$, then $I = 1$, $2$ or $4$. 
\item[(2)] If $G_{N} \simeq F_{384}$, then $I = 1$, $2$ or $4$. 
\end{list} 
\end{lemma}

\begin{proof} Let us show the assertion (1). We may assume that $X$ is 
projective. Since $Q_{16} < F_{128}$ by (2.6), one has 
$\text{rank}\, T(X) = 2$ by (4.4). 
Thus $I = 1$, $2$, $4$, $3$ or $6$ by (2.5). If $I = 6$, then 
$\mu_{3} < \mu_{I}$ and $H := \alpha^{-1}(\mu_{3})$ 
is a K3 group such that $H_{N} = F_{128}$ and $I = 3$. So, it suffices 
to show that 
$I \not= 3$. Assume that $I = 3$. Then, by (5.1), 
$G = F_{128} : \langle g \rangle$ 
where $\alpha(g) = \zeta_{3}$ and $\text{ord}(g) = 3$. Since 
$[F_{128}, F_{128}] \simeq C_{2} \times D_{8}$ by (2.6) and since the 
commutator subgroup 
is a characteristic subgroup, we have a new K3 group 
$K := (C_{2} \times D_{8}) : 
\langle g \rangle$ 
such that $K_{N} \simeq C_{2} \times D_{8}$ and $I = 3$. Let $c_{g}$ be the conjugate action 
by $g$ on $C_{2} \times D_{8}$. Since $C_{2} \times D_{8}$ has exactly one 
subgroup 
isomorphic to $C_{2} \times C_{4}$, we have a new K3 group 
$H = (C_{2} \times C_{4}) : \langle g \rangle$ such that $H_{N} = C_{2} \times C_{4}$ 
and $I = 3$. Since $C_{2} \times C_{4}$ contains exactly four order $4$ 
elements, $c_{g}$ 
fixes one of them, say $\tau$. Since there is then exactly 
two invulutions $\sigma$ such that  
$C_{2} \times C_{4} = \langle \sigma, \tau \rangle$, the conjugate 
action $c_{g}$ also fixes 
one of such $\sigma$. Hence, $H = (C_{2} \times C_{4}) \times \langle g \rangle$. Consider 
the element $h = \tau g$. Then $\text{ord}\, h = 12$ and 
$\alpha(h) = \zeta_{3}$. 
Let $M_{i}$ be the set of type (i) points of $X^{h}$ in (5.1)(3). Then, 
by (5.1)(3), one of $M_{i}$ 
is a one-point set, say $M = \{P\}$. Since $H$ is commutative, we have 
$a(P) = P$ 
for all $a \in H_{N}$. However, one would then have 
$$C_{2} \times C_{4} = H_{N} < \text{SL}(T_{X, P}) \simeq \text{SL}(2, \mathbf C)\,\,  ,$$ 
a contradiction to the fact that finite abelian subgroups of $\text{SL}(2, \mathbf C)$ 
must be cyclic.

Let us show the assertion (2). Note that $F_{384} = \langle F_{128}, 
\tau \rangle$ 
for some element $\tau$ of order $3$, and the Sylow $2$-subgroups of 
$F_{384}$ are 
exactly $F_{128}$, $\tau^{-1}F_{128}\tau$ and $\tau^{-2}F_{128}\tau^{2}$. For the same 
reason as in (1), it suffices to show that 
$I \not= 3$. Assume that $I = 3$. Then, by (5.1), 
$G = F_{384} : \langle g \rangle$ 
where $\alpha(g) = \zeta_{3}$ and $\text{ord}(g) = 3$. Consider a Sylow 
$3$-subgroup 
$H$ of $G$ containing $\tau$. Since $\vert G \vert = 2^{9}\cdot 3^{2}$, we have 
$\vert H \vert = 3^{2}$. 
Since $\vert G_{N} \vert = 2^{9}\cdot 3$, there is an element $h \in H$ 
such that $\alpha(h) = \zeta_{3}$. This element $h$ also acts by the 
conjugate on the set 
$\{F_{128}, \tau^{-1}F_{128}\tau , \tau^{-2}F_{128}\tau^{2}\}$ of 
Sylow $2$-subgroups of 
$G_{N}$. Thus, replacing $h$ by $h\tau^{i}$ if necessary, we have $h^{-1}F_{128}h = F_{128}$, 
and a new K3 group $K = \langle F_{128}, h \rangle$. 
Since $\alpha(h) = \zeta_{3}$ and $h$ is an element of a 
$3$-group $H$, 
we have $\text{ord}\, h = 3$ and 
$K = F_{128} : \langle h \rangle$ by (5.1)(1). 
However, this contradicts (5.2)(1). 
\end{proof} 
In the rest of this section, we prove (3.2)(2) and (4). 
\vskip0.5cm
\noindent
{\it Proof of (3.2)(2)}
\vskip0.3cm

The next proposition is obtained by [Mu] in the course of his proof of (2.3). 
For the notation of groups, we follow [Mu].   

\begin{proposition} \label{proposition:solv}[Mu, proposition 5.2 
and theorem 5.5] Let 
$G_{N}$ be a solvable symplectic K3 group. Then, $G_{N}$ and its order (indicated 
by [*]) is either  
\begin{list}{}{
\setlength{\leftmargin}{10pt}
\setlength{\labelwidth}{6pt}
}
\item[(I)] $2$-group $\,[\, 2^{n}\, ,\, 0 \leq n \leq 7\,]\,$; 
\item[(II)] $2 \cdot 3$-group $\,[\, 2^{n}3\, ,\, 0 \leq n \leq 7\,]\,$; 
moreover, if it is nilpotent, then $G_{N}$ is isomorphic to 
$C_{3}$, $C_{6}$ or 
$C_{2} \times C_{6}$; 
\item[(III)] $9\, \vert\, \vert G_{N} \vert$ and $G_{N}$ is either 
$$C_{3}^{2}\,[9]\,  ;\, A_{3,3}\, ,\, C_{3} \times S_{3}\, [18]\, ;
\, S_{3} \times S_{3}\, ,\, 
C_{3}^{2}:C_{4}\, ,\, A_{4} \times C_{3}\, [36]
\, ;\,N_{72}\, ,\, M_{9}\, ,\, A_{4,3}\, [72]\, ;$$ 
$$A_{4} \times A_{4}\, [144]\,  ,\, A_{4,4}\, [288]\, ;$$ 
\item[(IV)] $5\, \vert\, \vert G_{N} \vert$ and $G_{N}$ is either 
$$C_{5}\, [5]\, ,\, D_{10} (= C_{5} :C_{2})\, [10]\,  ,\, 
C_{5} : C_{4}\, [20]\, ,\, 
C_{2}^{4}:C_{5}\, [80]\, ,$$ 
$$C_{2}^{4} : D_{10}  [160]\,  ;$$ 
\item[(V)] $7 \vert \vert G_{N} \vert$ and $G_{N}$ is either 
$C_{7}\, [7]\, ,\,  C_{7} :C_{3}\, [21]$.  
\end{list} 
\end{proposition} 

Let us show (3.2)(2) dividing into the 5 cases (I) - (V) in (5.3). By 
(5.3), we may assume that $I \geq 2$. Then $X$ is projective as well.

First we consider the case where $G_{N}$ lies in (III), (IV), (V). 

\begin{lemma} \label{lemma:transval} 
\begin{list}{}{
\setlength{\leftmargin}{10pt}
\setlength{\labelwidth}{6pt}
}
\item[(1)] If $G_{N}$ is in the case (III), then $I \leq 12$. 
\item[(2)] If $G_{N}$ is in the case (IV),  then $I \leq 12$. 
\item[(3)] If $G_{N}$ is in the case (V),  then $I \leq 6$. 
\end{list} 
\end{lemma}

\begin{proof} First we shall show (1). Choose a subgroup $C_{3}^{2} \simeq 
\langle \tau_{1}, \tau_{2}  \rangle < G_{N}$. Then, by (4.5), one has 

$$ \text{rank}\, H^{2}(X, \mathbf Z)^{G_{N}} \leq  
\text{rank}\, H^{2}(X, \mathbf Z)^{\langle \tau_{1}, \tau_{2}  \rangle } 
= \frac{24 + 6 \times 8}{9} - 2 = 6\, .$$ 

Thus $\text{rank}\, T(X) \leq 5$ and we have $I \leq 12$ by (2.5). 

Proofs of (2) (resp. (3)) are similar if we choose a subgroup 
$C_{5}$ (resp. $C_{7}$)  in $G_{N}$: 

$$ \text{rank}\, H^{2}(X, \mathbf Z)^{G_{N}} \leq  
\text{rank}\, H^{2}(X, \mathbf Z)^{C_{5}} 
= \frac{24 + 4 \times 4}{5} - 2 = 6\,  ;$$  

$$ \text{rank}\, H^{2}(X, \mathbf Z)^{G_{N}} \leq  
\text{rank}\, H^{2}(X, \mathbf Z)^{C_{7}} 
= \frac{24 + 3 \times 6}{7} - 2 = 4\,  .$$  

\end{proof} 

Thus, when $G_{N}$ is in (III), (IV), or (V), we see 
that $\vert G \vert = \vert G_{N} \vert \cdot I < 2^{9}\cdot 3$ 
unless $G_{N}$ is one of 
$$(i)\, C_{2}^{4} : D_{10}\,  ,\, (ii)\, A_{4} \times A_{4}\,  ,\, (iii)\, A_{4,4}\, .$$ 

In case (i), we have $C_{2}^{4} : C_{5} \simeq H < G_{N}$. 
Here the order structure of 
$H$, which is also a subgroup (with no order $10$ element) 
of the affine transformation group 
$\mathbf F_{2}^{4} : \text{GL}(4, \mathbf F_{2})$, is as follows: 

$$
\begin{tabular}{*{4}{|c}{|}} 
\hline 
\text{order} & 1 & 2 & 5\\ 
\hline 
\text{cardinality}  & 1 & 15 & 64\\
\hline 
\end{tabular} 
$$

Then, one has 

$$ \text{rank}\, H^{2}(X, \mathbf Z)^{G_{N}} \leq  
\text{rank}\, H^{2}(X, \mathbf Z)^{H} 
= \frac{24 + 8 \cdot 15 + 4 \cdot 64}{80} - 2 = 3\,  .$$ 
Thus, $\text{rank}\, T(X) = 2$ and $I \leq 6$. 
Hence $\vert G \vert < 160 \cdot 6 = 960 < 2^{9} \cdot 3$.

Similarly, in case (ii), using the order structure of $G_{N} = A_{4} \times A_{4}$ below

$$
\begin{tabular}{*{5}{|c}{|}} 
\hline 
\text{order} & 1 & 2 & 3 & 6\\ 
\hline 
\text{cardinality}  & 1 & 15 & 80 & 48\\
\hline 
\end{tabular} 
$$

one can calculate 

$$ \text{rank}\, H^{2}(X, \mathbf Z)^{G_{N}} 
= \frac{24 + 8 \cdot 15 + 6 \cdot 80 + 2 \cdot 48}{144} - 2 = 3\,  .$$ 
Thus, $\text{rank}\, T(X) = 2$ and $I \leq 6$. 
Hence $\vert G \vert < 144 \cdot 6 = 864 < 2^{9} \cdot 3$.

Note that $A_{4} \times A_{4} < 
A_{4, 4} (:= (S_{4} \times S_{4}) \cap A_{8})$. Then, 
from the calculation above, we also find that $I = 1, 2, 3, 4$ or $6$ for 
$G_{N} = A_{4,4}$. 
Note that $\vert A_{4, 4} \vert \cdot 4 = 1152 < 2^{9} \cdot 3$, but 
$\vert A_{4, 4} \vert \cdot 6 = 
1728 > 2^{9} \cdot 3$. However, we can show the following: 

\begin{lemma} \label{lemma:excep} If $G_{N} \simeq A_{4, 4}$, then $I \not = 3$, $6$. 
\end{lemma}

\begin{proof} As in (5.2), it suffices to show that $I \not= 3$. Assume that $I = 3$. 
Then, by (5.1), $G = A_{4, 4} : \langle g \rangle$ 
where $\alpha(g) = \zeta_{3}$ and $\text{ord}(g) = 3$. Since 
$[A_{4, 4}, A_{4,4}] \simeq A_{4} \times A_{4}$ and since the commutator 
subgroup 
is a characteristic subgroup, we have a new K3 group $H := 
(A_{4} \times A_{4}) : 
\langle g \rangle$ 
such that $H_{N} \simeq A_{4} \times A_{4}$ and $I = 3$. 
Note that $A_{4} = C_{2}^{2} : C_{3}$ so that 
$H_{N} = A_{4} \times A_{4} = C_{2}^{4} : C_{3}^{2}$. 
Let $H_{3}$ be a Sylow $3$-subgroup 
of $H$ containing $C_{3}^{2}$. Since $\vert H \vert = 2^{4}\cdot 3^{3}$, 
we have 
$\vert H_{3} \vert = 3^{3}$. Note that $H_{3}$ acts on $H_{N}$ 
by the conjugate, 
say $\rho$. Since $C_{2}^{4}$ is the normal Sylow $2$-subgroup of $H_{N}$ 
(so that a characteristic subgroup of $H_{N}$), the conjugate action $\rho$ 
makes $C_{2}^{4}$ stable, and we have a group 
homomorphism 
$$\rho : H_{3} \longrightarrow \text{Aut}(C_{2}^{4}) \simeq \text{GL}(4, \mathbf F_{2}) \,\, .$$
Here $\vert \text{GL}(4, \mathbf F_{2}) \vert = 2^{6} \cdot 3^{2} \cdot 5 \cdot 7$. Thus, there 
is a non-trivial element $h \in \text{Ker}\, \rho$. Since $C_{3}^{2} (= H_{3} \cap H_{N})$ 
acts on $C_{2}^{4}$ faithfully, this $h$ satisfies $\alpha(h) = \zeta_{3}$ (after replacing $h$ by 
$h^{-1}$ if necessary). Moreover, $\text{ord}(h) = 3^{n}$ (by $h \in H_{3}$), 
and we have $\text{ord}(h) 
= 3$ by (5.1). Thus, we obtain a new K3 group $K = C_{2}^{4} \times \langle h \rangle$ 
such that $K_{N} = C_{2}^{4}$ and $I = 3$.  
Let $\sigma$ be an involution in $C_{2}^{4}$. Then $h \sigma$ is of order $6$ 
and satisfies $\alpha(h\sigma) = \zeta_{3}$. 
Let $M_{i}$ be the set of type (i) fixed points of $X^{h \sigma}$ 
described in (5.1)(2). Then, by (5.1)(2), one of $M_{i}$ 
is an at most $2$-point set, say $M = \{P, Q\}$. Since $K$ is commutative, 
we have $a(\{P, Q\}) = \{P, Q\}$ 
for all $a \in K_{N}$. Then, one would have an index $2$ subgroup 
$C_{2}^{3} < K_{N} (= C_{2}^{4})$ 
such that
$C_{2}^{3} < \text{SL}(T_{X, P}) \simeq \text{SL}(2, \mathbf C)$,  
a contradiction to the fact that finite abelian subgroups of 
$\text{SL}(2, \mathbf C)$ 
must be cyclic. \end{proof} 

Next, we consider the case (I), i.e. the case where $G_{N}$ is a $2$-group. 
Set 
$\vert G_{N} \vert = 2^{n}$. By (2.3) and (2.4), we have 
$G_{N} < F_{128}$ (as abstract groups). 
In particular, 
$n \leq 7$ and if $n = 7$, then $G_{N} \simeq F_{128}$. So,  by taking (5.2) 
into account, 
it suffices to show that $\vert G \vert < 2^{9}$ if $n \leq 6$.

Let us first consider the case where $G_{N}$ has an order $8$ element, 
say $\tau$. 
In this case, we have 
$$ \text{rank}\, H^{2}(X, \mathbf Z)^{G_{N}} \leq  
\text{rank}\, H^{2}(X, \mathbf Z)^{\langle \tau \rangle} 
= \frac{24 + 8 \times 1 + 4 \times 2 + 2 \times 4}{8} - 2 = 4\,  .$$ 

Thus, $I \leq 6$ and we have $\vert G \vert \leq 2^{6} \cdot 6 < 2^{9}$.

Next we consider the case where $G_{N}$ has no element of order $8$. Then, 
we have 
the following order structure of $G_{N}$:  

$$
\begin{tabular}{*{4}{|c}{|}} 
\hline 
\text{order} & 1 & 2 & 4\\ 
\hline 
\text{cardinality}  & 1 & $2k+1$ & $2m$\\
\hline 
\end{tabular} 
$$

where $k + m = 2^{n-1} - 1$. Moreover, $k \leq 17$ by $G_{N} < F_{128}$ and 
by (2.6)(1). 

If $n = 6$, then $k + m + 1 = 2^{5}$ and one has 
$$ \text{rank}\, H^{2}(X, \mathbf Z)^{G_{N}} =  
\frac{24 + 8(2k+1) + 4\cdot 2m}{2^{6}} - 2 = 2 + \frac{24 + 8k}{2^{6}} < 5\, 
,i.e.\, \leq  4\, .$$ 
Here the last inequality is because $k \leq 17$. Hence, 
$\text{rank}\, T(X) \leq 3$ 
and we have $I \leq 6$. Thus  
$\vert G \vert \leq 2^{6}\cdot 6 < 2^{9}$.

If $n = 5$, then $k + m = 15$ and $k \leq 15$. Thus, one has 
$$ \text{rank}\, H^{2}(X, \mathbf Z)^{G_{N}} =  
\frac{24 + 8(2k+1) + 4\cdot 2m}{2^{5}} - 2 = 
2 + \frac{24 + 8k}{2^{5}} < 7\, ,i.e.\,  \leq  6\, .$$ 
Hence, $\text{rank}\, T(X) \leq 5$ 
and we have $I \leq 12$. Thus  
$\vert G \vert \leq 2^{5}\cdot 12 < 2^{9}$.

Assume that $n \leq 4$. Then, if $\vert G \vert \geq 2^{9}$, we have $I \geq 2^{5} = 32$. 
In this case, one can check that $\varphi(I) \geq 12$ (see, for instance, 
the explicit list in [MO]). Then, $\vert G_{N} \vert 
\leq 2$ by the next lemma. We have then $\vert G \vert \leq 2 \cdot 66 < 2^{9}$. 

\begin{lemma} \label{lemma:large} Let $G$ be a K3 group on $X$. 
If $\varphi(I) \geq 12$. Then $\vert G_{N} \vert \leq 2$. 
\end{lemma} 

\begin{proof} By $\varphi(I) \geq 12$ and by (2.5), we have 
$\text{rank}\, T(X) \geq 12$. Let $g$ 
be a non-trivial element of $G_{N}$. 
Then $g^{*} \vert T(X) = id$ and $g$ fixes at least one ample class. 
Thus, 
$$\text{tr}\, (g^{*}\vert \text{NS}(X)) \geq 1 + (-1)\cdot 
(22 - \text{rank}\, T(X) -1) 
= \text{rank}\, T(X)  - 20 \,\, .$$

We also note that this inequality is strict if $\text{ord}\, g = 3$. 
Combining this with the topological Lefschetz fixed point formula, one has 
$$\vert X^{g} \vert = e(X^{g}) = 2 + \text{tr}\, (g^{*}\vert \text{NS}(X)) + 
\text{tr}\, (g^{*}\vert \text{T}(X)) \geq 
2\, \text{rank}\, T(X)  - 18 \geq 6\,\,  .$$ 

Thus, $g$ is an involution by (2.2) and by the remark above. 
Then, $G_{N} \simeq C_{2}^{n}$ for some $n$ and one has by (4.5)
$$ \text{rank}\, H^{2}(X, \mathbf Z)^{G_{N}} =  
\frac{24 + 8(2^{n}-1)}{2^{n}} - 2 = 6 + \frac{16}{2^{n}}\, .$$ 

Since $\text{rank}\, T(X) < \text{rank}\, H^{2}(X, \mathbf Z)^{G_{N}}$, we have then 
$$ 6 + \frac{16}{2^{n}} > 12\,\, ,i.e.\,\, n = 0\, ,\, 1\, .$$ 
\end{proof} 

Finally we consider the case (II), i.e. the case where $G_{N}$ is of order 
$2^{n} \cdot 3$. 
Then $n \leq 7$ and if $n = 7$, we have $G_{N} \simeq F_{384}$ by (2.3). 
The case $n = 7$ is settled by (5.2). Let $H$ be a Sylow $2$-subgroup of 
$G_{N}$. 
Then $\vert H \vert = 2^{n}$. 
By the argument in the case (I) and by the fact that 
$H^{2}(X, \mathbf Z)^{G_{N}} \subset H^{2}(X, \mathbf Z)^{H}$, we have 
$I \leq 6$ if $n = 6$ and $I \leq 12$ if $n = 5$. 
Then $\vert G \vert < 2^{9} \cdot 3$ for $n = 5$, $6$. 
Assume that $n \leq 4$.  Then, if 
$\vert G \vert (= \vert G_{N} \vert \cdot I) \geq 2^{9}\cdot 3$, 
then we have $I \geq 32$ 
and $\varphi (I) \geq 12$. Then, by (5.6), we would 
have $\vert G_{N} \vert \leq 2$, 
a contradiction. Thus $\vert G \vert < 2^{9} \cdot 3$ as well. 
Now we are done.   

\vskip0.5cm
\noindent
{\it Proof of (3.2)(4)}
\vskip0.3cm

Let us show the assertion (3.2)(4). Since $G$ is nilpotent, 
$G_{N}$ is also nilpotent. The previous argument for the solvable case already 
settled the case when $G_{N}$ is in (I) of (5.3). If a nilpotent group $G_{N}$ 
is in (II) of (5.3), then $\vert G_{N} \vert \leq 12$ and 
$$ \text{rank}\, H^{2}(X, \mathbf Z)^{G_{N}} \leq 
\text{rank}\, H^{2}(X, \mathbf Z)^{C_{3}}=  
\frac{24 + 6\cdot 2}{3} - 2 = 10\,\, .$$ 
Thus, $\text{rank}\, T(X) \leq 9$ 
and $\varphi(I) \leq 8$. This implies $I \leq 30$ (see eg. 
an explicit list in [MO]). We have then   
$\vert G \vert \leq 12\cdot 30 < 2^{9}$. 
If $G_{N}$ 
is in (III), (IV) or (V) of (5.3), then $G_{N}$ is either $C_{5}$, 
$C_{7}$ or $C_{3}^{2}$. 
(Recall that a nilpotent group must be the direct product 
of its Sylow subgroups.) 
Thus 
$\vert G_{N} \vert \leq 9$. Hence by (5.4), 
we have $\vert G \vert \leq 9 \cdot 12 < 2^{9}$ as well. 
Now we are done.

\section{Invariant polarization of a maximal nilpotent K3 group}
 
In this section, we shall prove proposition (3.3) along a similar 
line to [Ko2] and [OZ]. First,  we recall some 
basic facts on the Niemeier 
lattices needed in our arguments.
As in [OZ], our main reference concerning Niemeier lattices and 
their relations
with Mathieu groups is [CS, Chapters 10, 11, 16, 18]. 
\begin{definition} \label{definition:Neiemeier} The even negative 
definite 
unimodular lattices of rank 24
are called Niemeier lattices. 
There are exactly 24 isomorphism classes of the Niemeier lattices and
each isomorphism class is uniquely determined by its root lattice $R$, 
i.e. the sublattice generated by all the roots, the elements $x$ with
$x^{2} = -2$.
Except the Leech lattice which contains no root,
the other 23 lattices are the over-lattices of their root lattices. 
\end{definition}
We denote the Niemeier lattice $N$ whose root lattice $R$ by 
$N(R)$. Among 24 Niemeier lattices, the most relevant one for us is 
$N(A_{1}^{\oplus 24})$. Two other Niemeier lattices $N(A_{2}^{\oplus 12})$ 
and $N(A_{3}^{\oplus 8})$ will also appear in our argument.

Let $N = N(R)$ be a non-Leech Niemeier lattice. Denote by $O(N)$ (resp. by 
$O(R))$ the group of isometries of $N$ 
(resp. of $R$) and by $W(N) = W(R)$ the
Weyl group generated by the reflections given by the roots of $N$. 
Here $O(N) < O(R)$ and $W(N)$ is a normal subgroup of both 
$O(N)$ and $O(R)$. The invariant hyperplanes of
the reflections divide $N \otimes \mathbf R$ into finitely many chambers.
Each chamber is a fundamental domain of the action of $W(R)$. 
Fix a basis $\mathcal R := \{r_{i}\}_{i=1}^{24}$ of $R$ 
consisting of simple roots. The quotient group 
$S(N) := O(N)/W(R)$ is then identified with a subgroup of 
the full symmetry group $S(R) := O(R)/W(R)$ of the distinguished chamber 
${\mathcal C} :=
\{x \in N \otimes \mathbf R \vert (x, r) > 0\, ,\, r \in \mathcal R \}$, 
or a bit more concretely, $S(N)$ and $S(R)$ are subgroups of 
a larger group $S_{24}$ as: 
$$S(N) = \{g \in S(R)\, \vert g(N/R) = N/R \} < S(R) = \text{Aut}_{\text{graph}} (\mathcal R) < \text{Aut}_{\text{set}}(\mathcal R) = S_{24}\,\, ,$$
where the action by $S(R)$ on $N/R (\subset R^{*}/R)$ is induced by 
the natural action on $R^{*}/R$.  Here and hereafter, we denote by $M^{*}$ 
the dual lattice of a non-degenerate lattice $M$ and regard $M$ naturally 
as a submodule of finite index of $M^{*}$.

The groups $S(N)$ are explicitly calculated in [CS, Chapters 18, 16].
(See also [Ko1].) We need the following:
\begin{proposition} \label{proposition:root} [CS, Chapters 18, 16] 
Let $N$ be a 
non-Leech 
Niemeier lattice. Then, 
\begin{list}{}{
\setlength{\leftmargin}{10pt}
\setlength{\labelwidth}{6pt}
}
\item[(1)] $S(N) = M_{24}$ if
$N = N(A_1^{\oplus 24});$
\item[(2)] $S(N) = C_{2}.M_{12}$ if
$N = N(A_2^{\oplus 12});$
\item[(3)] $S(N) = C_{2}:(C_2^{\oplus 3} : L_3(2))$ if
$N = N(A_3^{\oplus 8});$ and 
\item[(4)] for other $N$, $S(N)$ is a subgroup of either 
$2.S_{6}$ or $3.S_{6}$. 
\end{list}
\end{proposition} 
Let us add a few remarks about the groups in (6.2)(1)-(3).

In case (1), i.e. the case where $N = N(R)$ and $R = A_{1}^{\oplus 24}$, 
we observe that 
$$\mathcal C_{24} := N/R \simeq \bold{F}_{2}^{\oplus 12} \subset 
R^{*}/R = \oplus_{i=1}^{24} \bold{F}_{2} 
\overline{r}_{i} \simeq \mathbf F_{2}^{\oplus 24}\, .$$ 
Here $\overline{r}_{i} := r_{i}/2 \, \text{mod}\, \mathbf Z r_{i}$. 
We note that $\mathcal R = \{r_{i}\}_{i=1}^{24}$ 
forms a Dynkin diagram of type 
$A_{1}^{\oplus 24}$.

Let $\mathcal P(\mathcal R)$ be the power set of $\mathcal R$. Then, we can 
identify $\mathcal P(\mathcal R)$ 
with $R^{*}/R$ by 
the following bijective correspondence: 
$$\iota : \mathcal P(\mathcal R) \ni A \mapsto \overline{r}_{A} := 
\frac{1}{2} \sum_{r_{j} \in A}r_{j} 
\, (\text{mod} \, R) 
\in R^{*}/R = (A_{1}^{\oplus 24})^{*}/A_{1}^{\oplus 24}.$$

In what follows, we freely identify these two sets, and 
we define $\vert x \vert$ ($x \in R^{*}/R$) 
to be the cardinality of $\iota^{-1}(x)$.

Then, under the identification 
by $\iota$, it is well known that $\emptyset, \mathcal R \in \mathcal C_{24}$ 
and that if $A \in \mathcal C_{24}$ ($A \not= \mathcal R$, $\emptyset$) 
then $\vert A \vert$ is either $8$, $12$, or $16$. We call 
$A \in \mathcal C_{24}$ an Octad (resp. a Dodecad) 
if $\vert A \vert = 8$ (resp. $12$). Note that $B \in \mathcal C_{24}$ with 
$\vert B \vert = 16$ is of the form $\mathcal R - A$ for some Octad $A$. 
It is also well known 
that the set of Octads forms a Steiner system $St(5, 8, 24)$ of $\mathcal R$ 
and generates $\mathcal C_{24}$ as $\mathbf F_{2}$-linear space. In this case, 
the embeddings $S(N) < S(R) < S_{24}$ explained above coincide with 
the natural inclusions 
$M_{24} < S_{24} = S_{24}$ for $N = N(A_{1}^{\oplus 24})$.

In the second case, the Mathieu group $M_{12} = W(N)/C_{2}$ acts naturally 
on the set of 12 
connected 
components of the Dynkin diagram 
$A_{2}^{\oplus 12}$ and $C_{2}$ interchanges the two vertices of all 
the components. We also 
note that $\vert M_{12} \vert = 2^{6} \cdot 3^{3} \cdot 5 \cdot 11$.

In the third case, we identify (non-canonically) the set of eight connected 
components of the Dynkin diagram
$A_{3}^{\oplus 8}$ with the three-dimensional linear space
$\mathbf F_{2}^{\oplus 3}$
over $\mathbf F_{2}$ by letting one connected component to be $0$. The group
$C_{2}:(C_{2}^{\oplus 3} : L_3(2))$ is the semi-direct product,
where $C_{2}$ interchanges the two edges of all the components,
$C_{2}^{\oplus 3}$ is the group of the parallel transformations of the affine
space
$\mathbf F_{2}^{\oplus 3}$ and $L_{3}(2) (\simeq L_{2}(7))$ is the linear
transformation group of $\mathbf F_{2}^{\oplus 3}$.

As in [Ko2] and [OZ], the next embedding theorem 
due to Kondo [Ko1] is an important ingredient in our proof: 

\begin{theorem} \label{theorem:kondo} [Ko1, Lemmas 5, 6] 
Let $K$ be a symplectic K3 group on $X$. Set $L := H^{2}(X, \mathbf Z)$, 
$L^{K} := \{x \in L\, \vert\, h(x) = x\, (\forall h \in K)\, \}$ and 
$L_{K} := \{y \in L\,  \vert\, (y, x) = 0\, (\forall x \in L^{K})\, \}$. 
Then, 

\begin{list}{}{
\setlength{\leftmargin}{10pt}
\setlength{\labelwidth}{6pt}
}
\item[(1)] There is a non-Leech Niemeier lattice $N$ such that 
$L_{K} \subset N$.
Moreover, the action by $K$ on $L_{K}$ extends to an action on $N$ so that
$L_{K} \simeq N_{K}$ and that $N^{K}$ contains a root, say $r^{0}$. 
Here the sublattices $N^{K}$ and $N_{K}$ of $N$ are defined in the same way 
as $L^{K}$ and $L_{K}$ of $L$.  
\item[(2)] Take $\mathcal R$ so that $r^{0} \in \mathcal R$. 
Then, the group action by $K$ on $N$ preserves 
the distinguished Weyl chamber 
$\mathcal C$ with respect to $\mathcal R$, and 
the naturally induced homomorphism $K \rightarrow S(N)$ is injective.
\end{list}
\end{theorem} 

\begin{corollary} \label{corollary:rank} [Ko1] Under the notation of 
(6.3), 
one has: 

\begin{list}{}{
\setlength{\leftmargin}{10pt}
\setlength{\labelwidth}{6pt}
}
\item[(1)] $\text{\rm rank}\, N^{K} = \text{\rm rank}\, L^{K} + 2$. 
\item[(2)] $(L^{K})^{*}/L^{K} \simeq (N^{K})^{*}/N^{K}$. In particular, 
$\vert \text{\rm det}\, N^{K} \vert = 
\vert \text{\rm det}\, L^{K} \vert$. 
\end{list}
\end{corollary} 

\begin{proof} The assertion (1) follows 
from $\text{rank}\, N^{K} = 24 - \text{rank}\, N_{K}$, 
$\text{rank}\, L^{K} = 22 - \text{rank}\, L_{K}$, and 
$N_{H} \simeq L_{H}$. Since $L$ and $N$ are unimodular and 
since the embeddings $L^{K} \subset L$ and $N^{K} \subset N$ are primitive, 
we have a natural isomorphisms $(L^{K})^{*}/L^{K} \simeq (L_{K})^{*}/L_{K}$ 
and $(N^{K})^{*}/N^{K} \simeq (N_{K})^{*}/N_{K}$. 
Now the result follows from $L_{K} \simeq N_{K}$. For the last equality, we 
may just note that $\vert \text{det}\, M \vert = \vert M^{*}/M \vert$ 
for a non-degenerate lattice $M$. \end{proof}  

We are now ready to prove (3.3). Let $G$ 
be a K3 group on $X$ such that $\vert G \vert = 2^{9}$. We denote by 
$K := G_{N}$, 
the symplectic part and by $I$ the transcendental value. By (3.2), 
$K = F_{128}$ (as abstract groups) 
and $I = 4$. In particular, $X$ is projective. We have also 
$\text{rank}\, L^{K} = 3$ by $K = F_{128} > Q_{16}$ and by (3.1)(1). 
Thus, $\text{rank}\, T(X) = 2$ and $\text{NS}(X)^{G} = 
\text{NS}(X)^{K} = \mathbf Z H$ for some ample class $H$. 
As in (6.3), we set $L := H^{2}(X, \mathbf Z)$. We shall fix these 
notations until the end of this section.

It remains to show $(H^{2}) = 4$. This will be completed in 
(6.11).

Let us first determine the Niemeier
lattice $N$ for our $K$. 

\begin{lemma} \label{lemma:ourG} The Niemeier lattice $N$ in (6.3) 
for
$K$ is $N(A_{1}^{\oplus 24})$.
\end{lemma}

\begin{proof} By (6.3)(2), $|S(N)|$ must be divied by 
$\vert K \vert = 2^{7}$.
Thus, $N$ is either $N(A_{1}^{\oplus 24})$, $N(A_{2}^{\oplus 12})$ or 
$N(A_{3}^{\oplus 8})$ 
by (6.2). Suppose that the second case occurs. Since $K$ fixes at least 
one element in $\mathcal R$ by (6.3)(1), we have $K < M_{12}$. However, 
this is impossible, because 
$\vert K \vert = 2^{7}$ and $\vert M_{12} \vert = 2^{6} \cdot k$ 
($(2, k) = 1$). Suppose that the third case occurs. Again for the same 
reason above, 
we have $K < C_{2}.L_{3}(2)$. However, this is impossible, because 
$\vert K \vert = 2^{7}$ and $\vert C_{2}.L_{3}(2) \vert = 2^{4} \cdot k'$ 
($(2, k') = 1$). Now we are done. 
\end{proof}

From now we set $N := N(A_{1}^{\oplus 24})$, $R := A_{1}^{\oplus 24}$ 
and take $\mathcal R = \{r_{i}\}_{i=1}^{24}$ as in (6.3)(2). 
By (6.2) and (6.3)(2) and (6.5), we have
$$K < M_{24} < S_{24} = \text{Aut}_{\text{graph}}(\mathcal R) 
= \text{Aut}_{\text{set}}(\mathcal R)\, .$$

\begin{lemma} \label{lemma:orbit} The orbit decomposition type of 
$\mathcal R$ 
by 
$K$ is $[1, 1, 2, 4, 16]$. 
\end{lemma}

\begin{proof} Note that $\text{rank}\, R^{K} = \text{rank}\, N^{K} =  5$ by 
$\text{rank}\, L^{K} = 3$ 
and by (6.4)(1). Thus $\mathcal R$ is divided into exactly 5 $K$-orbits. 
Since $K$ is a $2$-group and $K$ fixes at least one element by (6.2)(1), 
the orbit decomposition type 
is of the form $[1, 2^b, 2^c, 2^d, 2^e]$. We may assume that $0 \le b \le c 
\le d \le e$. In addition, $1 + 2^b + 2^c+ 2^d + 2^e = \vert \mathcal R \vert 
= 24$. It is now easy to see that $(b, c, d, e) = (0, 1, 2, 4)$. 
\end{proof}

By (6.6), after re-numbering of the elements of $\mathcal R$, we have 
$$R^{K} = \mathbf \langle s_{1}, s_{2}, s_{3}, s_{4}, s_{5} \rangle$$ 
where 
$$s_{1} = r_{1}\, ,\, s_{2} = r_{2}\, ,\, s_{3} = r_{3} + r_{4}\, ,\, 
s_{4} = r_{5} + \cdots + r_{8}\, ,\, 
s_{5} = r_{9} + \cdots +  r_{24}\, .$$ 

\begin{lemma} \label{lemma:inv}  
$$N^{K} = \mathbf \langle s_{1}, s_{2}, s_{3}, \frac{s_{1} + s_{2} + s_{3} + s_{4}}{2},  
\frac{s_{5}}{2} \rangle \,\,  .$$ In particular, 
$(N^{K})^{*}/N^{K} \simeq \mathbf Z/ 4 \oplus \mathbf Z /8 \oplus 
\mathbf Z /8$. 
\end{lemma}

\begin{proof} Since 
$R^{K} \subset N^{K} \subset (R^{*})^{K} = 
\mathbf \langle s_{1}/2\, ,\,  
s_{2}/2\, , \, s_{3}/2\, ,\, s_{4}/2\, ,\, s_{5}/2 \rangle$,   
the lattice $N^{K}$ is generated by $R^{K}$ and by the (representatives of) 
$K$-invariant elements of $\mathcal C_{24}$. Let us find out all such elements 
in $\mathcal C_{24}$. In what follows, we freely identify $\mathcal C_{24}$ 
with a subset of $\mathcal P(\mathcal R)$ by $\iota$, 
as it is explained after (6.2). 
By the shape of 
the orbit decomposition (6.6), there is no $K$-invariant Dodecad. Moreover, 
for the same reason, if there is a $K$-invariant Octad, then it must be 
$$(r_{1} + r_{2} + r_{3} + r_{4} + r_{5} + r_{6} + r_{7} + r_{8})/2  = 
(s_{1} + s_{2} + s_{3} + s_{4})/2\,\, ,i.e.$$ 
$$\{r_{1}\, ,\, r_{2}\, ,\, r_{3}\, ,\, r_{4}\, ,\, r_{5}\, ,\, r_{6}\, ,\, 
r_{7}\, ,\, r_{8}\}\,\, .$$ 
Let us show that this is indeed an Octad, i.e. 
an element of $\mathcal C_{24}$. 
Recall 
that the set of Octads of $\mathcal C_{24}$ forms a Steiner system 
$St(5, 8, 24)$ of $\mathcal R$. Then, there is an Octad 
$A \in \mathcal C_{24}$ 
containing $5$-element set 
$\{r_{1}, r_{5}, r_{6}, r_{7}, r_{8}\}$. Note that 
$$K(\{r_{1}, r_{5}, r_{6}, r_{7}, r_{8}\}) 
= \{r_{1}, r_{5}, r_{6}, r_{7}, r_{8}\}$$ 
by (6.6) and by $s_{1}, s_{4} \in R^{K}$. Then, by the Steiner 
property, we have 
$A = g(A)$ for all $g \in K$. Thus, this $A$ is a $K$-invarinat Octad. The argument here shows the existence of a $K$-invariant Octad. So, 
the only one possible candidate 
$(s_{1} + s_{2} + s_{3} + s_{4})/2$ is indeed a $K$-invariant Octad.
 
Since length-$16$ element of $\mathcal C_{24}$ is the complement of 
an Octad, it follows that  
$\{r_{9}, r_{10}, \cdots , r_{16}\}$ 
is the unique length-$16$, $K$-invariant element of $\mathcal C_{24}$. 
Hence  
$$N^{K} =  \mathbf \langle s_{1}, s_{2}, s_{3}, s_{4}, s_{5}\, \frac{s_{1} + s_{2} + s_{3} + s_{4}}{2},  
\frac{s_{5}}{2} \rangle\, ,$$
that is,  
$$N^{K} = \mathbf \langle s_{1}, s_{2}, s_{3}, \frac{s_{1} + 
s_{2} + s_{3} + s_{4}}{2},  
\frac{s_{5}}{2} \rangle\,\, .$$ 
The intersection matrix of $N^{K}$ with respect to this basis is

$$\begin{pmatrix} 
2&0&0&1&0\cr 
0&2&0&1&0\cr
0&0&4&2&0\cr
1&1&2&4&0\cr
0&0&0&0&8
\end{pmatrix}\,\, ,$$ 
and the elementary divisors of this matrix is $(1, 1, 4, 8, 8)$. 
This implies the result. 
\end{proof}  

\begin{lemma} \label{lemma:invII} 
\begin{list}{}{
\setlength{\leftmargin}{10pt}
\setlength{\labelwidth}{6pt}
}
\item[(1)] $(L^{K})^{*}/L^{K} \simeq \mathbf Z/ 4 \oplus \mathbf Z /8 \oplus \mathbf Z /8$. In particular, $\vert \text{det}\, L^{K} \vert = 2^{8}$. 
\item[(2)] If $x \in L^{K}$, then $(x^{2}) \equiv 0\, \text{mod}\, 4$. 
\end{list}
\end{lemma}
\begin{proof} The assertion (1) follows from (6.7) and (6.4)(2). 
By $\text{rank}\, L^{K} = 3$ and by (1), we can 
choose 
an integral basis $\langle f_{1}, f_{2}, f_{3} \rangle$ of $(L^{K})^{*}$ so that 
$\langle 4f_{1}, 8f_{2}, 8f_{3}\rangle$ forms an integral basis of $L^{K}$. 
For 
$x = x_{1}\cdot 4f_{1} + x_{2}\cdot 8f_{2} + x_{3}\cdot 8f_{3}$ ($x_{i} \in \mathbf Z$), one has 
$$(x^{2}) = 4x_{1}(x, f_{1}) + 8x_{2}(x, f_{2}) + 8x_{3}(x, f_{3}) \in 4 
\mathbf Z\, .$$
This implies the second assertion. 
\end{proof} 

\begin{lemma} \label{lemma:tpart} With respect to a suitable 
integral basis $\langle v_{1}, v_{2}\rangle$ of $T(X)$, 
the intersection matrix of $T(X)$ becomes the following form:
$$\begin{pmatrix} 
4m&0\cr 
0&4m
\end{pmatrix}\,\, \text{for\, some}\,\, m \in \mathbf Z\,\, .$$
\end{lemma}
\begin{proof} By (2.5), we have an isomorphism 
$T(X) \simeq \mathbf Z[\sqrt{-1}]$ as 
$\mathbf Z[\sqrt{-1}]$-modules. 
Since $\sqrt{-1}$ acts on the integral basis $\langle e_{1} :=  1, 
e_{2} := \sqrt{-1} \rangle$ of $\mathbf Z[\sqrt{-1}]$  
as $e_{1} \mapsto e_{2}$, $e_{2} \mapsto -e_{1}$, 
the group $G/K = \langle g\, 
\text{mod}\, K 
\rangle \simeq \mu_{4}$ 
acts on the corresponding integral basis $\langle v_{1}, v_{2} \rangle$ of $T(X)$ by $g^{*}(v_{1}) = v_{2}$ and $g^{*}(v_{2}) = -v_{1}$. 
Thus $(v_{1}, v_{2}) = (g^{*}(v_{1}), g^{*}(v_{2})) = (v_{2}, -v_{1})$, 
and $(v_{1}, v_{2}) = 0$. 
Similarly, $(v_{1}, v_{1}) = (g^{*}(v_{1}), g^{*}(v_{1})) = 
(v_{2}, v_{2})$. 
The result now follows 
from these two equalities and (6.8)(2). 
\end{proof}

\begin{lemma} \label{lemma:index} Set $l := [L^{K} : \mathbf Z H 
\oplus T(X)]$. 
Then $l = 1$ or $2$. Moreover, if $l = 2$, then 
$$ L^{K} = \mathbf Z \langle \frac{H + v_{1} + v_{2}}{2}, v_{1}, v_{2}\rangle\,\, .$$
Here $\langle v_{1}, v_{2} \rangle$ is an integral basis of $T(X)$ as in (6.9).
\end{lemma}

\begin{proof} Assume that $l \geq 2$. Since $T(X)$ and $H$ are both 
primitive in $L^{K}$, 
we have 
$$L^{K} = \mathbf Z \langle v_{1}, v_{2}, \frac{H + bv_{1} + cv_{2}}{l} 
\rangle$$ 
for some integers $b, c$ such that $0 \leq b, c \leq l-1$. By considering 
the action by 
$G/K = \langle g\, \text{mod}\, G_{N} 
\rangle \simeq \mu_{4}$, we have 

$$L^{K} \ni g^{*}(\frac{H + bv_{1} + cv_{2}}{l}) = 
\frac{H + bv_{2} - cv_{1}}{l}\,\, .$$ 
Here we used $g^{*}(v_{1}) = v_{2}$, $g^{*}(v_{2}) = -v_{1}$ (by the proof 
of (6.9)) and the fact that $H$ is $G$-invariant. Thus 

$$\frac{(b+c)v_{1} + (c - b) v_{2}}{l} \in L^{K}\,\, .$$ 
Since $T(X)$ is primitive in $L^{K}$, we have $l \vert b + c$ and $l \vert c-b$. Thus, 
$l \vert 2b$ and $l \vert 2c$. Since $0 \leq b, c \leq l -1$, it follows 
that $2b = 0$ or $2b = l$ 
and $2c = 0$ or $2c = l$. If $b = c = 0$, then $H/l \in L^{K}$, 
a contradiction to the 
primitivity of $H$. If $2b = l$ and $c = 0$, then 
$H/l + v_{1}/2 \in L^{K}$. Applying $g$, 
we have $H/l + v_{2}/2 \in L^{K}$. Thus $(v_{1} - v_{2})/2 \in L^{K}$, 
a contradiction 
to the primitivity of $T(X) \subset L^{K}$. We get a similar 
contradiction when $b = 0$ and 
$2c = l$. Hence $b = c = l/2$ and we have $H/l + (v_{1} + v_{2})/2 \in L^{K}$. 
Since $v_{1} + v_{2} \in L^{K}$, it follows that $2H/l \in L^{K}$. 
Thus $l = 2$ (when $l \geq 2$) by the primitivity of $H$, 
and $b = c = 1$. 
\end{proof} 

The next lemma completes the proof of (3.3). 
  
\begin{lemma} \label{lemma:pol} $(H^{2}) = 4$. 
\end{lemma}
 
\begin{proof} By (6.8)(2), we can write $(H^{2}) = 4n$ for some positive 
integer $n$. We need to show that $n=1$. Let $m$ be a positive integer 
in (6.9). 

First consider the case where $l= 2$ (Here $l$ is the index defined in 
(6.10)). In this case, we have by (6.8)(1)  
$$4\cdot 2^{8} = l^{2} \cdot \text{det}\, L^{K} = (H^{2}) 
\cdot \text{det}\, T(X) 
= 4n \cdot 16m^{2}\, .$$ 
Thus $nm^{2} = 16$. Moreover, by $(H + v_{1} + v_{2})/2 \in L^{K}$ and by (6.8)(2), we have 
$$n + 2m = ((H + v_{1} + v_{2})/2)^{2}) \equiv 0\, \text{mod}\, 4\, .$$ 
Thus $(m, n) = (2, 4)$ and the intersection matrix of $L^{K}$ (with respect 
to the basis in (6.10)) becomes 
$$\begin{pmatrix}
8&4&4\cr 
4&8&0\cr
4&0&8
\end{pmatrix}\,\, .$$
However, the elementary divisors of this matrix is $(4, 4, 16)$, 
a contradiction to (6.8)(1). 

Next consider the case where $l = 1$. In this case, we have  
$$2^{8} = \text{det}\, L^{G_{N}} = (H^{2}) \cdot \text{det}\, T(X) = 
4n \cdot 16m^{2}\, .$$ 
Thus $nm^{2} = 4$ and $(m, n)$ is either $(1, 4)$ or $(2, 1)$. 
Assume that $(m, n) = (1, 4)$. Then, the 
intersection 
matrix of $L^{K} = \mathbf Z H \oplus T(X)$ (with respect to the basis 
$\langle H, v_{1}, v_{2} \rangle$) becomes 
 
$$\begin{pmatrix}
16&0&0\cr 
0&4&0\cr
0&0&4
\end{pmatrix}\,\, .$$

However, the elementary divisors of this matrix is $(4, 4, 16)$, 
a contradiction to (6.8)(1). Thus $(m, n) = (2, 1)$ and we are done.  
\end{proof} 

\end{document}